\documentclass[11pt, reqno]{amsart}
\usepackage{amsfonts,latexsym}
\usepackage{amsmath}
\usepackage{amscd}
\usepackage{float,amsmath,amssymb,mathrsfs,bm,multirow,graphics}
\usepackage[dvips]{graphicx}
\usepackage[percent]{overpic}
\usepackage[pdftex]{color}

\newcommand{\nequation}{\setcounter{equation}{0}}
\renewcommand{\theequation}{\mbox{\arabic{section}.\arabic{equation}}}

\newcommand{\C}{{\Bbb C}}

\newcommand{\proofbegin}{\noindent{\it Proof.\,\,}}
\newcommand{\proofend}{\hfill$\Box$\bigskip}

\allowdisplaybreaks

\newtheorem{theorem}{Theorem}[section]
\newtheorem{proposition}[theorem]{Proposition}
\newtheorem{lemma}[theorem]{Lemma}

\newtheorem{remark}[theorem]{Remark}

\input epsf
\title[The Unified Method II]{\sc The Unified Method: II~NLS on the Half-Line with $t$-Periodic Boundary Conditions}
\author{J. Lenells}
\address{J.L.: Department of Mathematics, Baylor University, One Bear Place \#97328, Waco, TX 76798, USA.}
\email{Jonatan\_Lenells@baylor.edu}

\author{A. S. Fokas}
\address{A.S.F.: Department of Applied Mathematics and Theoretical Physics, University of Cambridge, Cambridge CB3 0WA, UK, and Research Center of Mathematics, Academy of Athens, 11527, Greece.}
\email{T.Fokas@damtp.cam.ac.uk}

\begin{document}
\begin{abstract} 
\noindent Boundary value problems for integrable nonlinear evolution PDEs formulated on the half-line can be analyzed by the unified method introduced by one of the authors and used extensively in the literature. The implementation of this general method to this particular class of problems yields the solution in terms of the unique solution of a matrix Riemann-Hilbert problem formulated in the complex $k$-plane (the Fourier plane), which has a jump matrix with explicit $(x,t)$-dependence involving four scalar functions of $k$, called spectral functions. Two of these functions depend on the initial data, whereas the other two depend on all boundary values. The most difficult step of the new method is the characterization of the latter two spectral functions in terms of the given initial and boundary data, i.e. the elimination of the unknown boundary values. For certain boundary conditions, called linearizable, this can be achieved simply using algebraic manipulations. Here, we first present an effective characterization of the spectral functions in terms of the given initial and boundary data for the general case of non-linearizable boundary conditions. This characterization is based on the analysis of the so-called global relation and on the introduction of the so-called Gelfand-Levitan-Marchenko representations of the eigenfunctions defining the spectral functions.
We then concentrate on the physically significant case of $t$-periodic Dirichlet boundary data. 
After presenting certain heuristic arguments which suggest that the Neumann boundary values become periodic as $t\to\infty$, we show that for the case of the NLS with a sine-wave as Dirichlet data, the asymptotics of the Neumann boundary values can be computed explicitly at least up to third order in a perturbative expansion and indeed at least up to this order are asymptotically periodic. 
\end{abstract}

\maketitle



\tableofcontents

\section{Introduction}\nequation
This is the second in a series of papers addressing the most difficult problem in the analysis of integrable nonlinear evolution PDEs, namely the problem of expressing the so-called {\it spectral functions} in terms of the given initial and boundary conditions. In \cite{FLnonlinearizable} this problem was analyzed {\it directly}, i.e. without employing the so-called Gelfand-Levitan-Marchenko (GLM) representations. Here, we revisit the same problem, but we employ the GLM representations; it appears that this latter approach might have an advantage for studying problems with $t$-periodic boundary problems. Furthermore, for the NLS we also present certain explicit formulas for the large $t$-behavior of the Neumann boundary values corresponding to the simple case of a sine-wave as Dirichlet data.

We refer the interested reader to \cite{FLnonlinearizable} for an introduction to the unified method of \cite{F1997, F2002} and for a discussion of the difference between linearizable versus nonlinearizable boundary value problems. Here we only emphasize that the unified method expresses the solution $q(x,t)$ of an integrable nonlinear evolution PDE formulated on the half-line as an integral in the complex $k$-plane. This representation is similar to the representation obtained by the new method for the linearized version of the given PDE, but also involves the entries of a certain matrix-valued function $M(x,t,k)$, which is the solution of a matrix Riemann-Hilbert (RH) problem. The main advantage of the unified method is that this RH problem involves a jump matrix with {\it explicit} $(x,t)$-dependence, uniquely defined in terms of four scalar functions called spectral functions and denoted by $\{a(k), b(k), A(k), B(k)\}$. The functions $\{a(k), b(k)\}$ are defined in terms of the initial data $q_0(x) = q(x,0)$, whereas the functions $\{A(k), B(k)\}$ are defined in terms of the boundary values. For example, for the nonlinear Schr\"odinger (NLS) equation
\begin{align}\label{nls}
 & i\frac{\partial q}{\partial t} + \frac{\partial^2 q}{\partial x^2} - 2 \lambda |q|^2 q = 0, \qquad \lambda = \pm 1,
\end{align}
the functions $\{A(k), B(k)\}$ are defined in terms of $g_0(t) = q(0,t)$ and $g_1(t) = q_x(0,t)$. However, for a well-posed boundary value problem only one of the boundary values (or their combination) is prescribed as a boundary condition; for example, for the Dirichlet problem of the NLS, $g_0(t)$ is given but $g_1(t)$ is unknown. The problem of determining $g_1(t)$ in terms of $q_0(x)$ and $g_0(t)$ is known as the problem of determining the Dirichlet to Neumann map. 

For boundary conditions which decay for large $t$, by utilizing the crucial feature of the new method that the dependence of the jump matrix on $(x,t)$ is explicit, it is possible to obtain useful asymptotic information about the solution {\it without} the need to characterize the spectral functions $\{A(k), B(k)\}$ in terms of the given initial and boundary conditions (or equivalently without the need of constructing the Dirichlet to Neumann map). However, the complete solution of non-linearizable boundary value problems requires the characterization of $A(k)$ and $B(k)$ in terms of the given initial and boundary conditions. The effective solution of this problem is particularly important for the physically significant case of boundary conditions which are periodic in $t$, since in this case it is not possible to obtain the rigorous form of the long time asymptotics, without the full characterization of $A(k)$ and $B(k)$ (in spite of this difficulty important results about the large $t$-asymptotics of such problems are presented in \cite{BIK2007, BIK2009, BKSZ2010}).

The construction of the Dirichlet to Neumann map can be achieved by analyzing the so-called {\it global relation}, which is an algebraic equation coupling $\{a(k), b(k), A(k), B(k)\}$.
A breakthrough in the analysis of the global relation was achieved in \cite{BFS2003}, where the global relation was solved {\it explicitly}  for the unknown Neumann boundary values in terms of the given initial conditions, the given Dirichlet boundary conditions, and the functions $\{L_j(t,s), M_j(t,s)\}_{j=1}^2$ appearing in the Gelfand-Levitan-Marchenko (GLM) representations of the eigenfunctions $\{\Phi_1(t,k), \Phi_2(t,k)\}$. The numerical implementation of an improved version of these formulae was presented in \cite{Z2006}, see also \cite{Z2007}. Following these important developments, a similar analysis for the sine-Gordon, the two versions of the mKdV and the two versions of the KdV, was carried out in \cite{F2005} and \cite{TF2008}. 

\subsection{A new formulation of the Dirichlet to Neumann map}
Here we first present an extension of the above results. In particular, for the NLS we show that the global relation actually yields {\it two linear equations coupling} $\{L_j(t,s)\}_1^2$ and $\{M_j(t,s)\}_1^2$. The equation for the Neumann boundary values obtained in \cite{BFS2003} can be obtained by evaluating one of these two equations at $s = t$. Similar results are valid for other integrable nonlinear PDEs like the sG, mKdV, and KdV equations. 

By employing these two new linear equations, it is possible to obtain a simpler characterization of the generalized Dirichlet to Neumann map than the one presented in \cite{BFS2003, F2005,TF2008}. 

\subsection{Time-periodic boundary conditions}
Boundary conditions which are periodic in time are of great importance in applications. Examples include the  periodic motion of a wave maker or the vertical motion of a buoy in the ocean.
We will consider the following question: Assuming that the Dirichlet boundary conditions are periodic in time, do the Neumann boundary values approach a function that is periodic in time as $t \to \infty$? After presenting a heuristic argument why this question might be positively answered in the general case, we provide more concrete evidence by showing that: (a) For any Dirichlet data $g_0(t)$ such that $\dot{g}_0(t+t_p) - \dot{g}_0(t)$ decays faster than $1/\sqrt{t}$ as $t \to \infty$, $g_1(t)$ is asymptotically periodic at least to {\it first} order in the perturbative expansion.
(b) For the prototypical example when $g_0(t)$ is a sine-wave, the asymptotics of $g_1(t)$ to {\it third} order in perturbation theory can be determined explicitly and at least up to this order $g_1$ is asymptotically periodic.

\subsection{Organization of the paper}
In section \ref{nlssec}, we derive the new characterization of the Dirichlet to Neumann map for the NLS equation. Time-periodic boundary conditions are considered in section \ref{tperiodicsec}.

\section{NLS revisited}\label{nlssec}\nequation
Let the vector $(\Phi_1, \Phi_2)$ satisfy the `normalized' $t$-part of the usual Lax pair of the NLS evaluated at $x = 0$ (see \cite{Fbook} for details):
\begin{align}\nonumber
&\begin{pmatrix} \Phi_1(t,k) \\ \Phi_2(t,k) \end{pmatrix}_t + 4ik^2 \begin{pmatrix} 1 & 0 \\ 0 & 0 \end{pmatrix} \begin{pmatrix} \Phi_1(t,k) \\ \Phi_2(t,k) \end{pmatrix}
	\\ \nonumber
&\hspace{2cm} =  \left\{2k\begin{pmatrix} 0 & g_0(t) \\ \lambda \bar{g}_0(t) & 0 \end{pmatrix}
- i\lambda \begin{pmatrix} |g_0(t)|^2 & -\lambda g_1(t) \\
\bar{g}_1(t) & -|g_0(t)|^2 \end{pmatrix} \right\}
\begin{pmatrix} \Phi_1(t,k) \\ \Phi_2(t,k) \end{pmatrix},
	\\ \label{normalizedtpart}
& \begin{pmatrix}\Phi_1(0,k) \\ \Phi_2(0,k) \end{pmatrix} = \begin{pmatrix} 0 \\ 1 \end{pmatrix}, \qquad \qquad \qquad \qquad \qquad
 0 < t < T, \; k \in \C.
 \end{align}
 Then, it is shown in \cite{BK2000} that $(\Phi_1, \Phi_2)$ admits the GLM representations
 \begin{subequations}
 \begin{align}
\label{Phi1}
  & \Phi_1(t,k) = \int_{-t}^t \left[ L_1(t,s) - \frac{i}{2}g_0(t) M_2(t,s) + kM_1(t,s)\right] e^{2ik^2(s-t)} ds,
   	\\
 & \Phi_2(t,k) = 1 + \int_{-t}^t \left[ L_2(t,s) + \frac{i\lambda}{2}\bar{g}_0(t)M_1(t,s) + kM_2(t,s)\right] e^{2ik^2(s-t)} ds, 
  	\\ \nonumber
&  \hspace{6cm}0 < t < T, \; -t < s < t, \; k \in \C,
 \end{align}
 \end{subequations}
if and only if the functions $\{L_j(t,s), M_j(t,s)\}_1^2$ satisfy the following well-posed Goursat hyperbolic problem:
\begin{subequations}\label{Goursat}
  \begin{align}\label{Goursata}
  &  L_{1t} - L_{1s} = ig_1(t)L_2 + \alpha(t)M_1 + \beta(t) M_2,
    	\\ \label{Goursatb}
  &  M_{1t} - M_{1s} = 2g_0(t) L_2 + ig_1(t) M_2,
    	\\  \label{Goursatc}
  &  L_{2t} + L_{2s} = -i\lambda\bar{g}_1(t)L_1 - \alpha(t)M_2 + \lambda \bar{\beta}(t) M_1,
    	\\ \label{Goursatd}
  &  M_{2t} + M_{2s} = 2\lambda \bar{g}_0(t) L_1 - i\lambda \bar{g}_1(t) M_1, \qquad
    0 < t < T, \; -t < s < t,	
   \end{align}
 \end{subequations}
 \vspace{-.5cm}
\begin{align}\label{LMinitial}
  L_1(t,t) = \frac{i}{2}g_1(t), \quad M_1(t,t) = g_0(t), \quad L_2(t,-t) = M_2(t,-t) = 0,
\end{align}
where the functions $\alpha(t)$ and $\beta(t)$ are defined in terms of $g_0(t)$ and $g_1(t)$ by
\begin{align}\label{alphabetadef}
  \alpha(t) = \frac{\lambda}{2}\left(g_0(t) \bar{g}_1(t) - \bar{g}_0(t) g_1(t)\right), \qquad
  \beta(t) = \frac{1}{2}\left(i\dot{g}_0(t) - \lambda|g_0(t)|^2 g_0(t)\right).
\end{align}
The general solution of the equation
\begin{align}  \label{Fequation1}
  F_t(t,s) - F_s(t,s) = f(t,s), \qquad 0 < t < T, \; -t < s < t,
\end{align}
is given by (see appendix \ref{characteristicapp})
\begin{align}\label{Ffromcharacteristics1}
F(t,s) = F\left(\frac{t + s}{2}, \frac{t + s}{2}\right)
+ \int_{\frac{t+s}{2}}^tf(\tau, t + s- \tau) d\tau.
\end{align}
Thus, equations (\ref{Goursata}) and (\ref{Goursatb}) imply that we require the knowledge of the functions $L_1(x,x)$ and $M_1(x,x)$, $x = \frac{t+s}{2}$, which are given by the first two equations in (\ref{LMinitial}). Similarly, the general solution of the equation
\begin{align}
  F_t(t,s) + F_s(t,s) = f(t,s), \qquad 0 < t < T, \; -t < s < t,
\end{align}
is given by
\begin{align}\label{Ffromcharacteristics2}
  F(t,s) = F\left(\frac{t-s}{2}, -\frac{t-s}{2}\right)+ \int_{\frac{t-s}{2}}^t f(\tau, \tau + s -t)d\tau.
\end{align}
Thus, equations (\ref{Goursatc}) and (\ref{Goursatd}) imply that we require the knowledge of the functions $L_2(x,-x)$ and $M_2(x,-x)$, $ x = (t-s)/2$, which are given by the last two equations in (\ref{LMinitial}).

In the case of zero initial conditions, $s(k) = I$, and then the global relation (2.16) of section 2 of \cite{FLnonlinearizable} yields
\begin{align} \label{PhijFj}
  \Phi_1(t,k) = F_1(t,k), \quad  \Phi_2(t,k) = F_2(t,k), \qquad  0 < t < T, \; k \in \C^+,
\end{align}
where $F_1$ and $F_2$ are bounded and analytic functions of $k \in \C^+$ and of $O(1/k)$ as $k \to \infty$ in $\C^+$.

We replace in the first of equations (\ref{PhijFj}) $\Phi_1$ by the rhs of (\ref{Phi1}), we replace $s$ by $\tau$, we multiply the resulting equation by 
$$k e^{2ik^2(t-s)}, \qquad  -t < s < t,$$
and integrate with respect to $dk$ along $\partial D_1$ which is the oriented boundary of the first quadrant of the complex $k$-plane:
\begin{equation*}
\int_{\partial D_1} k \int_{-t}^{t} \left[ L_1(t,\tau)-\frac{i}{2}g_0(t)M_2(t,\tau)+kM_1(t,\tau) \right] e^{2ik^2(\tau-s)}d\tau dk
\end{equation*}
\begin{equation}\label{2.11}
=\int_{\partial D_1} ke^{2ik^2(t-s)}F_1(t,k) dk = 0.         
\end{equation}
We simplify this equation by using the following facts:
(a) $F_1(t,k)$ is analytic in $k$ in $\C^+$ and of $O(1/k)$ as $k \to \infty$ in $\C^+$. Thus, Jordan's lemma implies that the rhs of (\ref{2.11}) vanishes.
(b) By employing the change of variables $2k^2 = l$, we find
$$\int_{\partial D_1} ke^{2ik^2(\tau -s)}dk = \frac{\pi}{2}\delta(\tau - s),$$
thus the first two terms on the rhs of (\ref{2.11}) yield
$$\frac{\pi}{2}\left[L_1(t,s) - \frac{i}{2}g_0(t) M_2(t,s)\right].$$
(c) Integrating by parts the term involving $k^2$, we find
\begin{align}\label{intM1identity}
 & \int_{\partial D_1} k^2 e^{-2ik^2s}\left(\int_{-t}^t e^{2ik^2\tau}M_1(t, \tau)d\tau\right) dk
  	\\\nonumber
& \hspace{2cm} = \frac{1}{2i}\int_{\partial D_1}\left[e^{2ik^2(t-s)}M_1(t,t) - e^{-2ik^2(t+s)}M_1(t,-t)\right]dk
	\\ \nonumber
 & \hspace{2.5cm} - \frac{1}{2i}\int_{\partial D_1} \left(\int_{-t}^t e^{-2ik^2(s-\tau)}M_{1\tau}(t,\tau)d\tau\right) dk.
\end{align}
The exponential $\exp[2ik^2(t-s)]$ is bounded and analytic in $D_1$, thus Jordan's lemma implies that the first integral on the rhs of (\ref{intM1identity}) vanishes. Similarly, the integral with respect to $d\tau$ from $s$ to $t$ vanishes. Hence, the rhs of (\ref{intM1identity}) equals
$$\frac{ic_0}{2}\left[\frac{M_1(t,-t)}{\sqrt{s+t}} + \int_{-t}^s \frac{M_{1\tau}(t,\tau)}{\sqrt{s -\tau}} d\tau\right], \qquad c_0 = \int_{\partial D_1} e^{-2il^2} dl.$$
By deforming the contour $\partial D_1$ to the real axis we find that
\begin{align*}  
  c_0 = \sqrt{\frac{\pi}{2}}e^{-\frac{i\pi}{4}}.
\end{align*}
Hence, equation (\ref{2.11}) yields
\begin{subequations}
\begin{align}\label{preL1L2eqsa}
  L_1(t,s) - \frac{i}{2}g_0(t) M_2(t,s) + \frac{e^{\frac{i\pi}{4}}}{\sqrt{2\pi}} \left[\frac{M_1(t,-t)}{\sqrt{s+t}} + \int_{-t}^s \frac{M_{1\tau}(t,\tau)}{\sqrt{s-\tau}} d\tau\right] = 0.
  \end{align}
Similarly, the second of equations (\ref{PhijFj}) yields
\begin{align}\label{preL1L2eqsb}
  L_2(t,s) + \frac{i\lambda}{2} \bar{g}_0(t) M_1(t,s) + \frac{e^{\frac{i\pi}{4}}}{\sqrt{2\pi}} \int_{-t}^s \frac{M_{2\tau}(t, \tau)}{\sqrt{s-\tau}} d\tau = 0,
\end{align}
\end{subequations}
where we have used that $M_2(t,-t) = 0$.

Evaluating equation (\ref{preL1L2eqsb}) at $s = -t$ and using that $L_2(t,-t) = 0$, it follows that $M_1(t,-t) = 0$. Then, evaluating (\ref{preL1L2eqsa}) at $s = -t$, we find that $L_1(t,-t) = 0$.

In summary, we have derived the following result: 
\begin{proposition}
Let $\{L_j(t, s), M_j(t,s)\}_1^2$, $0 < t < T$, $-t < s < t$, be the GLM representations associated with the $t$-part of the Lax pair of the NLS equation (\ref{nls}) evaluated at $x = 0$, see equations (\ref{normalizedtpart})-(\ref{alphabetadef}). Then, the global relation of the NLS yields the following equations:
\begin{subequations}\label{L1L2eqs}
  \begin{align}\label{L1L2eqsa}
 &   L_1(t,s) - \frac{i}{2}g_0(t) M_2(t,s) + \frac{e^{\frac{i\pi}{4}}}{\sqrt{2\pi}} \int_{-t}^s \frac{M_{1\tau}(t,\tau) d\tau}{\sqrt{s - \tau}} = 0,
    	\\
 &   L_2(t,s) + \frac{i\lambda}{2}\bar{g}_0(t) M_1(t,s) + \frac{e^{\frac{i\pi}{4}}}{\sqrt{2\pi}} \int_{-t}^s \frac{M_{2\tau}(t,\tau) d\tau}{\sqrt{s - \tau}} = 0,\label{L1L2eqsb}
    	\\
& M_1(t,-t) = L_1(t,-t) = 0.\label{L1L2eqsc}	
    \end{align}
\end{subequations}
\end{proposition}

\begin{remark}\upshape
Evaluating equation (\ref{L1L2eqsa}) at $s = t$, we recover the basic formula of \cite{BFS2003, Fbook}:
\begin{align}\label{g1basic}
  g_1(t) = g_0(t) M_2(t,t) - \sqrt{\frac{2}{\pi}} e^{-\frac{i\pi}{4}} \int_{-t}^t \frac{M_{1\tau}(t, \tau) d\tau}{\sqrt{t - \tau}}.
\end{align}
\end{remark}

\begin{remark}\upshape
Equations (\ref{L1L2eqs}) encode the information contained in the global relation. The advantage of the GLM representations is that the symmetry relations associated with the global relation are automatically incorporated into this information. Indeed, letting $k \to -k$ yields nothing new at the level of GLM representations.
\end{remark}

For a function $F(t,s)$, $-t < s < t$, we denote by $\mathcal{A}$ the Abel transform in the second variable:
$$(\mathcal{A}F)(t,s) = \int_{-t}^s F(t,\tau) \frac{d\tau}{\sqrt{s - \tau}}.$$
The inverse $\mathcal{A}^{-1}$ of this transform is given by
$$(\mathcal{A}^{-1}F)(t,s) = \frac{1}{\pi} \frac{d}{ds}\int_{-t}^s \frac{F(t,\tau)d\tau}{\sqrt{s - \tau}} = \frac{1}{\pi}\left[\int_{-t}^s \frac{\frac{\partial F}{\partial \tau}(t,\tau)d\tau}{\sqrt{s-\tau}} + \frac{F(t,-t)}{\sqrt{t+s}}\right].$$

\begin{theorem}\label{th3}
Let $q(x,t)$ satisfy the NLS (\ref{nls}) on the half-line and assume that $q_0(x) = 0$, $x \geq 0$. In this case $M_1(t,-t)=M_2(t,-t)=0$, thus the inverse Abel transform of $\{M_j\}_1^2$ is given by
\begin{equation}
\mathcal{A} ^{-1} M_j(t,s)=\frac{1}{\pi} \int_{-t}^{s} \frac{\frac{\partial M_j}{\partial \tau}(t,\tau)d \tau}{\sqrt{s-\tau}}, \quad j=1,2, \ -t<s<t, \ 0<t<T.
\end{equation}
The function $g_1(t):=q_x(0,t)$ is given by
\begin{equation}\label{g1th3}
  g_1(t) = g_0(t) M_2(t,t) - \sqrt{2\pi} e^{-i\pi/4}(\mathcal{A}^{-1}M_1)(t,t) 
\end{equation}
where the functions $\{M_j(t,s)\}_1^2$, $|s| < t < T$, satisfy the following system of quadratically nonlinear integral equations:
\begin{subequations}\label{Msystem}
\begin{align} \nonumber
  M_{1}(t,s) = \;& g_0\Bigl(\frac{t +s}{2}\Bigr) + \int_{\frac{t+s}{2}}^t \biggl[- g_0(\tau)\left(i\lambda\bar{g}_0(\tau) M_1 +  \sqrt{2\pi} e^{i\pi/4} (\mathcal{A}^{-1}M_2) \right)
 	\\ \label{Msystema}
 & + i\left(g_0(\tau) M_2(\tau,\tau) - \sqrt{2\pi} e^{-i\pi/4}(\mathcal{A}^{-1}M_1)(\tau,\tau) \right) M_2 \biggr] d\tau,
 	 \\ \nonumber
 M_{2}(t,s) = \;& \lambda \int_{\frac{t-s}{2}}^t \biggl[ \bar{g}_0(\tau) \left(ig_0(\tau) M_2 - \sqrt{2\pi} e^{i\pi/4} (\mathcal{A}^{-1}M_1) \right) 
	\\ \label{Msystemb}
& - i \left(
\bar{g}_0(\tau) \bar{M}_2(\tau,\tau) - \sqrt{2\pi} e^{i\pi/4}(\mathcal{A}^{-1}\bar{M}_1)(\tau,\tau) \right) M_1\biggr] d\tau.
 \end{align} 
 \end{subequations}
Unless otherwise specified, the functions in the integrand of (\ref{Msystema}) are evaluated at $(\tau, t+s-\tau)$ and the functions in the integrand of (\ref{Msystemb}) are evaluated at $(\tau, \tau - t + s)$.
\end{theorem}
\proofbegin
The expression (\ref{g1th3}) for $g_1(t)$ follows from (\ref{g1basic}). The first two equations in (\ref{L1L2eqs}) can be written as
\begin{align}\nonumber
& L_1(t,s) = \frac{i}{2}g_0(t) M_2(t,s) - \sqrt{\frac{\pi}{2}} e^{i\pi/4} (\mathcal{A}^{-1}M_1)(t,s),
	\\ \nonumber
& L_2(t,s) = - \frac{i\lambda}{2}\bar{g}_0(t) M_1(t,s) -  \sqrt{\frac{\pi}{2}} e^{i\pi/4} (\mathcal{A}^{-1}M_2)(t,s).
 \end{align}
Using these equations together with (\ref{g1th3}) to eliminate $L_1$, $L_2$, and $g_1$ from equations (\ref{Goursatb}) and (\ref{Goursatd}), we find the following system for $M_1$ and $M_2$:
\begin{align*}
  (\partial_t - \partial_s)M_{1}(t,s)  = & - 2g_0(t)\left(\frac{i\lambda}{2}\bar{g}_0(t) M_1(t,s) +  \sqrt{\frac{\pi}{2}} e^{i\pi/4} (\mathcal{A}^{-1}M_2)(t,s) \right)
 	\\
 & + i\left(g_0(t) M_2(t,t) - \sqrt{2\pi} e^{-i\pi/4}(\mathcal{A}^{-1}M_1)(t,t) \right) M_2(t,s),
 	 \\
 (\partial_t + \partial_s)M_{2}(t,s) = \;& 2\lambda \bar{g}_0(t) \left(\frac{i}{2}g_0(t) M_2(t,s) - \sqrt{\frac{\pi}{2}} e^{i\pi/4} (\mathcal{A}^{-1}M_1)(t,s) \right) 
	\\
& - i \lambda \left(
\bar{g}_0(t) \bar{M}_2(t,t) - \sqrt{2\pi} e^{i\pi/4}(\mathcal{A}^{-1}\bar{M}_1)(t,t) \right) M_1(t,s).
 \end{align*}
The identities (\ref{Ffromcharacteristics1}) and (\ref{Ffromcharacteristics2}) together with the initial conditions (\ref{LMinitial}) yield (\ref{Msystem}).

\proofend


\subsection{An effective perturbative scheme}
In what follows we show that equations (\ref{g1th3})-(\ref{Msystem}) provide an {\it effective} characterization of the Dirichlet to Neumann map, namely we show that they yield a well-defined perturbative scheme for computing $g_1(t)$.

Substituting into (\ref{Msystem}) the expansions
\begin{align*}
& M_1 = \epsilon M_1^{(1)} + \epsilon^2 M_1^{(2)} + \cdots, \qquad
M_2 = \epsilon M_2^{(1)} + \epsilon^2 M_2^{(2)} + \cdots, 
	\\
& g_0 = \epsilon g_0^{(1)} + \epsilon^2 g_0^{(2)} + \cdots,
\end{align*}
where $\epsilon> 0$ is a small parameter, we find the following: The terms of $O(\epsilon)$ yield
\begin{align} \nonumber
  M_{1}^{(1)}(t,s) = g_0^{(1)}\Bigl(\frac{t +s}{2}\Bigr),
 \qquad M_{2}^{(1)}(t,s) = 0.
 \end{align} 
The terms of $O(\epsilon^2)$ yield
 \begin{align} \nonumber
 & M_{1}^{(2)}(t,s) = g_0^{(2)}\Bigl(\frac{t +s}{2}\Bigr),
 	 \\ \nonumber
 &M_{2}^{(2)}(t,s) = \lambda \sqrt{2\pi} e^{i\pi/4}  \int_{\frac{t-s}{2}}^t \biggl[-\bar{g}_0^{(1)}(\tau)  (\mathcal{A}^{-1}M_1^{(1)})(\tau, \tau - t + s)
 	\\ \nonumber
& \hspace{5cm} + i (\mathcal{A}^{-1}\bar{M}_1^{(1)})(\tau,\tau) M_1^{(1)}(\tau, \tau - t + s)\biggr] d\tau.
 \end{align} 
In general, the terms of $O(\epsilon^n)$ yield
 \begin{align} \nonumber
 & M_{1}^{(n)}(t,s) = g_0^{(n)}\Bigl(\frac{t +s}{2}\Bigr) + \text{known lower order terms},
 	 \\ \nonumber
 &M_{2}^{(n)}(t,s) = \text{known lower order terms}.
 \end{align} 
Since $g_0(t)$ is known, it follows that the system (\ref{Msystem}) can be solved perturbatively to all orders. Thus, for `small' boundary conditions, the system yields an effective perturbative characterization of $M_1$ and $M_2$, in which each term can be computed uniquely in a well-defined recursive manner.

\section{Time-periodic boundary conditions}\label{tperiodicsec}\nequation
In applications it is common to have boundary conditions which are periodic or at least asymptotically periodic in time. For the Dirichlet problem this means that $g_0(t + t_p) = g_0(t)$ for some period $t_p > 0$, or, in the asymptotically periodic case, that $g_0(t + t_p) - g_0(t)$ decays to zero as $t \to \infty$.
Intuitively we expect the effect of the initial data in the boundary region to diminish for large $t$, so that the leading asymptotics of the solution near the boundary should be determined by $g_0(t)$. This suggests that if $g_0(t)$ is periodic, then the solution $q(x,t)$ will eventually become periodic in time for small $x$ as $t \to \infty$. In this section, we will investigate the large time asymptotics of the Neumann data for the NLS equation, assuming that the Dirichlet data are periodic. In particular, we will consider the question whether (asymptotically) periodic Dirichlet boundary conditions lead to asymptotically periodic Neumann boundary values, i.e. whether the Dirichlet to Neumann map takes (asymptotically) periodic functions to asymptotically periodic functions. 
After presenting a heuristic argument why this might be true in the general case, we will provide more concrete evidence in the case of vanishing initial data by showing that: (a) For any Dirichlet data $g_0(t)$ such that $\dot{g}_0(t+t_p) - \dot{g}_0(t)$ decays faster than $1/\sqrt{t}$ as $t \to \infty$, $g_1(t)$ is asymptotically periodic at least to first order in perturbation theory.
(b) For the special case when $g_0(t)$ is a sine-wave, the asymptotics of $g_1(t)$ to third order in perturbation theory can be determined explicitly and at least up to third order $g_1$ is asymptotically periodic. 

\subsection{Heuristic argument}
In order to establish the asymptotic periodicity of $g_1(t)$, we must analyze the analogous property of the solution $\{M_1(t,s), M_2(t,s)\}$ of equations (\ref{Msystem}). These equations  involve the following three types of integrals:
\begin{align}\label{I1Mdef}
  & (I_1M)(t, t+s) := \int_{\frac{t+s}{2}}^t g(\tau)M(\tau, t+s-\tau)d\tau,
  	\\\label{I2Mdef}
  & (I_2M)(t, t-s) := \int_{\frac{t-s}{2}}^t g(\tau) M(\tau, s - t + \tau) d\tau,
  	\\\label{I3Mdef}
  & (I_3M)(t,s) := \int_{-t}^s \frac{M_\tau(t, \tau)}{\sqrt{s - \tau}} d\tau, \qquad 0 < t < T, \; -t < s < t.		
\end{align}
We assume that $g(t)$ is $2\pi$ periodic and analyze the properties of the integrals (\ref{I1Mdef})-(\ref{I3Mdef}), under the transformations
\begin{align}\label{tsshifts}
  t \to t + 2\pi, \qquad s \to s + 2\pi.
\end{align}

We start with the integral in (\ref{I1Mdef}):
\begin{align}
  (I_1M)(t + 2\pi, t + s + 4\pi) = \int_{\frac{t+s}{2} + 2\pi}^{t + 2\pi} g(\tau) M(\tau, t +s + 4\pi - \tau)d\tau.
\end{align}
Replacing $\tau$ with $\tau + 2\pi$ in the rhs of the above equation we find
\begin{align}
  (I_1M)(t + 2\pi, t + s + 4\pi) = (I_1\hat{M})(t, t+s),
\end{align}
where
\begin{align}
  \hat{M}(t,s) = M(t + 2\pi, s + 2\pi).
\end{align}

Similarly, replacing $\tau$ with $\tau + 2\pi$ in the integral resulting from (\ref{I2Mdef}) under the transformations (\ref{tsshifts}), we find
\begin{align*}
(I_2M)(t + 2\pi, t - s) =\;& \int_{\frac{t - s}{2} - 2\pi}^t g(\tau) M(\tau + 2\pi, s - t + \tau + 2\pi) d\tau
	\\
=\;& \int_{\frac{t-s}{2}}^t g(\tau) M(\tau + 2\pi, s - t + \tau + 2\pi)d\tau
	\\
& + \int_{\frac{t-s}{2} - 2\pi}^{\frac{t-s}{2}} g(\tau) M(\tau + 2\pi, s - t + \tau + 2\pi) d\tau.
\end{align*}
Replacing in the second integral in the rhs of the above equation $\tau$ with $(t-s)/2 + \tau$, we find
\begin{align}\label{I2Mshifted}
  (I_2M)(t + 2\pi, t&-s) = (I_2\hat{M})(t , t-s) 
  	\\ \nonumber
&  + \int_{-2\pi}^0 g\Bigl(\frac{t-s}{2} + \tau\Bigr)M\Bigl(\frac{t-s}{2} + \tau + 2\pi, \frac{s-t}{2} + \tau + 2\pi\Bigr) d\tau.
\end{align}
We note that for large $t$ the function $M$ appearing in the integral of the second term of (\ref{I2Mshifted}) can be approximated by $M(t, -t)$ and also recall the equations
$$L_1(t,-t) = M_1(t,-t) = M_2(t, -t) = 0.$$
Finally, we consider the integral in (\ref{I3Mdef}). Replacing $\tau$ with $\tau + 2\pi$ in the integral resulting from (\ref{I3Mdef}) under the transformations (\ref{tsshifts}), we find
\begin{align*}
 (I_3M)(t + 2\pi, s + 2\pi) & = \int_{-t - 4\pi}^s \frac{M_\tau(t + 2\pi, \tau + 2\pi)}{\sqrt{s - \tau}} d\tau
 	\\
& = \int_{-t}^s \frac{M_\tau(t + 2\pi, \tau + 2\pi)}{\sqrt{s - \tau}}d\tau
+ \int_{-t - 4\pi}^{-t} \frac{M_\tau(t + 2\pi, \tau + 2\pi)}{\sqrt{s - \tau}} d\tau.
\end{align*}
Replacing in the second integral in the rhs of the above equation $\tau$ with $\tau - t$, we find
\begin{align}\label{I3Mshifted}
  (I_3M)(t + 2\pi, s + 2\pi)
  = (I_3\hat{M})(t,s) + \int_{-4\pi}^0 \frac{M_\tau(t + 2\pi, -t + \tau + 2\pi)}{\sqrt{s + t -\tau}}d\tau.
\end{align}
We note that if $M_\tau$ is bounded in $t$, the integral in the rhs of (\ref{I3Mshifted}) vanishes as $t \to \infty$.

The above properties suggest that it might be possible to prove directly from equations (\ref{g1th3})-(\ref{Msystem}) that the periodicity of $g_0(t)$ implies the asymptotic periodicity of $g_1(t)$. However, the rigorous proof of this result remains open.

\subsection{Time-periodic boundary conditions in the linear limit}
\begin{theorem}\label{g11th}
  Let $q(x,t) = \epsilon q_1(x,t) + O(\epsilon^2)$ be the perturbative solution of the NLS equation (\ref{nls}) on the half-line with $q(x,0) = 0$, $x > 0$, and $q(0,t) = \epsilon g_{01}(t) + O(\epsilon^2)$, where $g_{01}(t)$ is a sufficiently smooth function compatible with the zero initial data, i.e. $g_{01}$ and its derivatives vanish at $t = 0$.
Let $g_1(t) = q_x(0,t) = \epsilon g_{11}(t) + O(\epsilon^2)$, denote the Neumann boundary values. If 
\begin{align}\label{g01asymptoticallyperiodic}
\dot{g}_{01}(t + t_p) - \dot{g}_{01}(t) = O\biggl(\frac{1}{t^{\frac{1}{2} + \nu}}\biggr), \qquad t \to \infty,
\end{align}
for some $t_p > 0$ and $\nu > 0$, then
\begin{align}\label{g11asymptoticallyperiodic}
  g_{11}(t + t_p) - g_{11}(t) = \begin{cases} O(t^{-\nu}), & 0 < \nu  < 1/2,\\
O(t^{-1/2}\ln{t}), & \nu = 1/2, \\
   O(t^{-1/2}), & \nu > 1/2, 
   \end{cases} \qquad t \to \infty.
\end{align}  
\end{theorem}
\proofbegin
The Neumann boundary values $g_{11}(t)$ are given by an Abel transform of $g_{01}(t)$, see \cite{FLnonlinearizable}:
\begin{align} \label{g11formula}
g_{11}(t) = - \frac{e^{-\frac{i\pi}{4}}}{\sqrt{\pi}} h(t), \quad \text{where}\quad h(t) := \int_0^t \frac{\dot{g}_{01}(s)}{\sqrt{t - s}} ds.
\end{align}
We write
\begin{align*}
h(t+t_p) = \int_0^{t+t_p} \frac{\dot{g}_{01}(s)ds}{\sqrt{t+t_p - s}}
= \int_0^{t_p} \frac{\dot{g}_{01}(s)ds}{\sqrt{t+t_p - s}}
+\int_{t_p}^{t+t_p} \frac{\dot{g}_{01}(s)ds}{\sqrt{t+t_p - s}}.
\end{align*}
Changing variables $s \to s + t_p$ in the second integral, we find
\begin{align*}
h(t+t_p) &= \int_0^{t_p}\frac{\dot{g}_{01}(s)ds}{\sqrt{t+t_p - s}} + \int_0^{t} \frac{\dot{g}_{01}(s + t_p)ds}{\sqrt{t - s}}
	\\
&= \int_0^{2\pi}\frac{\dot{g}_{01}(s)ds}{\sqrt{t+t_p - s}}
+ \int_0^{t} \frac{\dot{g}_{01}(s +t_p) - \dot{g}_{01}(s)}{\sqrt{t - s}}ds + h(t).
\end{align*}
Since the first integral on the rhs is of $O(t^{-1/2})$ as $t \to \infty$, we infer that
$$h(t+t_p) - h(t) = \int_0^{t} \frac{\dot{g}_{01}(s +t_p) - \dot{g}_{01}(s)}{\sqrt{t - s}}ds + O(t^{-1/2}), \qquad t \to \infty.$$
Suppose now that $g_{01}(t)$ satisfies (\ref{g01asymptoticallyperiodic}). Then, there exist constants $K > 0$ and $M > 0$, such that
\begin{align*}
\biggl|\int_0^t \frac{\dot{g}_{01}(s +t_p) - \dot{g}_{01}(s)}{\sqrt{t - s}} ds\biggr|
\leq \int_0^K \frac{|\dot{g}_{01}(s +t_p) - \dot{g}_{01}(s)|}{\sqrt{t - s}} ds 
+  M \int_K^t \frac{1}{s^{\frac{1}{2}+\nu} \sqrt{t - s}} ds
\end{align*}
for $t > K$. The first integral on the rhs is clearly of $O(t^{-1/2})$ as $t\to \infty$, whereas the second integral satisfies
\begin{align*}
\int_K^t \frac{ds}{s^{\frac{1}{2}+\nu} \sqrt{t - s}} 
& = 
   {}_2F_1\biggl(\frac{1}{2},\nu +\frac{1}{2}; \frac{3}{2}; 1-\frac{K}{t}\biggr)\frac{2\sqrt{1- K/t} }{t^{\nu}} 
   	\\
& = \begin{cases} O(t^{-\nu}), & 0 < \nu  < 1/2,\\
O(t^{-1/2}\ln{t}), & \nu = 1/2, \\
   O(t^{-1/2}), & \nu > 1/2,
   \end{cases} \qquad t \to \infty.
\end{align*}
This proves (\ref{g11asymptoticallyperiodic}).
\proofend

\subsection{A sine-wave as boundary data}
The prototypical example of a periodic boundary profile which vanishes at $t = 0$ is the sine-wave $g_0(t) = \sin{t}$.
For this example we can take the results of the previous subsection one step further: We will compute the asymptotics of $g_{1}(t)$ to {\it third} order in perturbation theory and, as expected, find that the result is asymptotically periodic.

\begin{theorem}\label{g13th}
Let 
$$q(x,t) = \epsilon q_1(x,t) + \epsilon^2 q_2(x,t) + \cdots, \qquad \epsilon \to 0,$$
be the perturbative solution of NLS equation (\ref{nls}) on the half-line with the initial conditions $q(x,0) = 0$, $x > 0$, and the Dirichlet boundary conditions 
\begin{align}\label{Dirichletsineexpansion}
  q(0,t) = \epsilon \sin{t} + O(\epsilon^4).
\end{align}  
Let 
$$g_1(t) = q_x(0,t) = \epsilon g_{11}(t) + \epsilon^2 g_{12}(t) + \epsilon^3 g_{13}(t) + O(\epsilon^4),$$
denote the corresponding Neumann boundary values.
Then 
\begin{subequations}\label{g11g12g13}
\begin{align}\label{g11g12g13a}
& g_{11}(t) = -\frac{e^{-\frac{i \pi }{4}}}{\sqrt{2}}(\cos {t} + \sin{t}) + O(t^{-3/2}), && t \to \infty,
	\\\label{g11g12g13b}
& g_{12}(t) = 0, && t > 0,
 	\\\label{g11g12g13c}
& \lambda g_{13}(t) = c_1 \cos{t} + c_2 \sin{t} + c_3(\cos{3t} - \sin{3t}) + O(t^{-1/16}), && t\to \infty,
\end{align}
\end{subequations}
where the constant $c_3$ is given by
$$c_3 = \frac{(1-\sqrt{3}) (\sqrt{3}+i)}{16},$$
whereas expressions for the constants $c_1$ and $c_2$ are given in equation (\ref{c1c2def}) below. 

Equations (\ref{g11g12g13}) imply that the Neumann boundary values are asymptotically periodic at least to third order in the perturbative expansion. More precisely,
\begin{subequations}\label{Deltag11g12g13}
\begin{align}
&g_{11}(t + 2\pi) - g_{11}(t) = O(t^{-3/2}), && t \to \infty,
	\\
&g_{12}(t + 2\pi) - g_{12}(t) = 0, && t > 0
\end{align}
and
\begin{align}
g_{13}(t + 2\pi) - g_{13}(t) = O(t^{-1/16}), && t \to \infty.
\end{align}
\end{subequations}
\end{theorem}

\begin{remark}\upshape
If higher order terms are included in (\ref{Dirichletsineexpansion}), i.e.
\begin{align*}
  q(0,t) = \epsilon \sin{t} + \epsilon^2 g_{02}(t) + \epsilon^3 g_{03}(t) + O(\epsilon^4),
\end{align*}  
where $g_{02}(t)$ and $g_{03}(t)$ are sufficiently smooth functions compatible with the zero initial data, then the only effect of the functions $g_{02}(t)$ and $g_{03}(t)$ is to add the terms
$$- \frac{e^{-\frac{i \pi }{4}}}{\sqrt{\pi }} \int_0^t \frac{\dot{g}_{02}(s)}{\sqrt{t - s}}ds \quad \text{and}\quad
 - \frac{e^{-\frac{i \pi }{4}}}{\sqrt{\pi }} \int_0^t \frac{\dot{g}_{03}(s)}{\sqrt{t - s}}ds$$
to the expressions for $g_{12}(t)$ and $g_{13}(t)$, respectively. Employing the argument used to prove theorem \ref{g11th}, we deduce that $g_1(t)$ is still asymptotically periodic to third order in the perturbative expansion provided that $g_{02}(t)$ and $g_{03}(t)$ are asymptotically periodic. 
\end{remark}

\noindent{\it Proof of theorem \ref{g13th}.\,\,}
Since $\lambda$ only affects the sign of $g_{13}$, we may assume that $\lambda =1$.
When $g_{01}(t) = \sin{t}$, the function $h(t)$ defined in (\ref{g11formula}), i.e. 
$$h(t) = \int_0^t \frac{\cos{s}}{\sqrt{t-s}} ds,$$
satisfies
\begin{align}\label{hasymptotics}
h(t) & = 2 \sqrt{t} \, _1F_2\left(1;\frac{3}{4},\frac{5}{4};-\frac{t^2}{4}\right)
	\\ \nonumber
&= \sqrt{\frac{\pi }{2}} (\sin{t} +\cos{t})-\frac{1}{2 t^{3/2}}+ \frac{15}{8 t^{7/2}}-\frac{945}{32 t^{11/2}} + O(t^{-15/2}), \qquad t\to \infty.
\end{align}
Since $g_{11} = - e^{-\frac{i\pi}{4}} h/\sqrt{\pi}$, this proves (\ref{g11g12g13a}). Since $g_{01}(t)$ has no term of $O(\epsilon^2)$, $g_{12}(t)$ vanishes identically, which proves (\ref{g11g12g13b}). Clearly, equations (\ref{Deltag11g12g13}) are a consequence of (\ref{g11g12g13}). Therefore, it only remains to prove (\ref{g11g12g13c}), and this relies on the following expression for $g_{13}(t)$, which is obtained by setting $g_{01}(t) = \sin{t}$ in the general expression for $g_{13}(t)$ derived in \cite{FLnonlinearizable}:
\begin{align} \nonumber
\lambda g_{13}(t) =\; & \frac{2 c}{\pi i}\int_0^t  \frac{\sin^3{t'}}{\sqrt{t - t'}} dt' 
-  \frac{2 c}{\pi i}\int_0^t \sin^2{t'} \int_0^{t'}\frac{\cos{t''}}{\sqrt{t - t''}}dt'' dt'
	\\\nonumber
& - \sqrt{\frac{2}{\pi}} \sin{t}\int_0^t h(t') \sin{t'}dt'
	\\ \label{g13formula}
& 
- \frac{ c}{\pi i} I_{ssc}(t)
 +  \frac{ c}{\pi^2} I_{shh}(t)
+ \frac{ c i}{\pi^2} I_{hsh}(t)
+ \frac{ c}{\pi^2} I_{hhs}(t)
	\\\nonumber
 =: &\; T_1 + T_2 + T_3 + T_4 + T_5 + T_6 + T_7,
\end{align}	
where $c = \frac{\sqrt{\pi}}{2}e^{-\frac{i\pi}{4}}$ and the functions $I_{ssc}$, $I_{shh}$, $I_{hsh}$, and $I_{hhs}$ are defined by
\begin{subequations}\label{Ixxxdef}
\begin{align} \label{Isscdef}
& I_{ssc}(t) = \int_0^t dt' \sin{t'} \int_0^{t'} dt'' \sin{t''} \int_0^{t''} dt''' \frac{\cos{t'''}}{(t-t'+t'' - t''')^{3/2}},
 	\\ \label{Ishhdef}
& I_{shh}(t) = \int_0^t dt' \sin{t'} \int_0^{t'} dt'' h(t'') \int_0^{t''}dt''' \frac{h(t''')}{(t - t' + t'' - t''')^{3/2}},
	\\ \label{Ihshdef}
& I_{hsh}(t) = \int_0^t dt' h(t') \int_0^{t'} dt'' \sin{t''} \int_0^{t''} dt''' \frac{h(t''')}{(t - t' + t'' - t''')^{3/2}},
	\\ \label{Ihhsdef}
& I_{hhs}(t) = \int_0^t dt' h(t') \int_0^{t'} dt'' h(t'') \int_0^{t''} dt'''\frac{\sin(t''')}{(t - t' + t'' - t''')^{3/2}}.
\end{align}
\end{subequations}
We will treat the different expressions in (\ref{Ixxxdef}) in a series of lemmas. 
In this respect, we introduce the following notations:
\begin{itemize}
\item $S(z)$ and $C(z)$ are the Fresnel integrals defined by
\begin{align}
S(z) = \int_0^z \sin\frac{\pi s^2}{2} ds, \qquad
C(z) = \int_0^z \cos\frac{\pi s^2}{2} ds.
\end{align}

\item Given a function $f(t)$, $Af(t)$ denotes the Abel transform,
\begin{align}\label{Adef}
Af(t) := \int_0^t \frac{f(s)ds}{\sqrt{t-s}}.
\end{align}

\item Given a function $f(t)$, $If(t)$ denotes the integral
\begin{align}
  If(t) := \int_0^t f(s) ds.
\end{align}
 
\item The functions $H_s(t)$ and $H_c(t)$ are defined by
\begin{align}\label{HcHsdef}
  H_s(t):= \int_0^{t} h(s) \sin{s} ds, \qquad H_c(t) := \int_0^{t} h(s) \cos{s} ds.
\end{align}  
\end{itemize}
Note that the inversion of the Abel transform yields $Ah(t) = \pi \sin{t}$. The function $h(t)$ can be written as
\begin{align}\label{hintermsofCS}
h(t) = \sqrt{2 \pi } \biggl(C\biggl(\sqrt{\frac{2}{\pi }} \sqrt{t}\biggr) \cos
   (t)+S\biggl(\sqrt{\frac{2}{\pi }} \sqrt{t}\biggr) \sin (t)\biggr).
\end{align}

\begin{lemma}\label{T1lemma}
The term $T_1$, i.e. the first term of the rhs of (\ref{g13formula}), satisfies
\begin{align*}
 T_1
& = \frac{i c}{3 \sqrt{2 \pi }}
\biggl\{-9 C\biggl(\sqrt{\frac{2}{\pi }} \sqrt{t}\biggr) \sin
   (t)+\sqrt{3} C\biggl(\sqrt{\frac{6}{\pi }} \sqrt{t}\biggr) \sin (3 t)
  	\\
  &\hspace{2cm} +9S\biggl(\sqrt{\frac{2}{\pi }} \sqrt{t}\biggr) \cos (t)-\sqrt{3}
   S\biggl(\sqrt{\frac{6}{\pi }} \sqrt{t}\biggr) \cos (3 t)\biggr\}.
\end{align*}
\end{lemma}
\proofbegin
\begin{align*}
  \int_0^t  \frac{\sin^3{t'}}{\sqrt{t - t'}} dt' 
  = \frac{1}{4} \int_0^t  \frac{3\sin{t'} - \sin{3 t'}}{\sqrt{t - t'}} dt'
  = \frac{1}{4} \int_0^t  \frac{3\sin(t-s) - \sin(3t-3s)}{\sqrt{s}} ds.  
\end{align*}
The result follows by expanding the sine functions and expressing the result in terms of $C(z)$ and $S(z)$.
\proofend

\begin{lemma}
The terms $T_2$ and $T_3$, i.e. the second and third terms of the rhs of (\ref{g13formula}),  satisfy
\begin{align*}
T_2 =\;& -\frac{2 c}{\pi i}\biggl[\biggl(\frac{t}{2} - \frac{\sin{2t}}{4}\biggr)h - A\biggl(\cos(t)\biggl(\frac{t}{2} - \frac{\sin{2t}}{4}\biggr)\biggr)\biggr]
\end{align*}
and $T_3 =  - \sqrt{\frac{2}{\pi}} H_s(t)\sin{t}$.
\end{lemma}
\proofbegin
This follows using integration by parts and the identity
$$\int_0^t \sin^2{s} \, ds = \frac{t}{2}-\frac{1}{4} \sin {2t}.$$
\proofend

In order to find the asymptotics of $T_4$-$T_7$, our strategy is to reduce the expressions for $I_{ssc}$, $I_{shh}$, $I_{hsh}$, $I_{hhs}$ to expressions involving only the Fresnel integrals $C$, $S$, the functions $H_c$, $H_s$, $h$, and the operators $A$ and $I$.

\begin{lemma}
The function $I_{ssc}(t)$ defined in (\ref{Isscdef}) can be expressed in the following form:
\begin{align*}
 I_{ssc} = \;& \frac{1}{8} \biggl\{\sqrt{6 \pi } \biggl(S\biggl(\sqrt{\frac{6}{\pi }}
   \sqrt{t}\biggr) \cos (3 t)-C\biggl(\sqrt{\frac{6}{\pi }}
   \sqrt{t}\biggr) \sin (3 t)\biggr)
   	\\
&   +\sqrt{2 \pi }
   C\biggl(\sqrt{\frac{2}{\pi }} \sqrt{t}\biggr) (3 \sin {t}+\sin (3 t)-2
   t \cos {t})
   	\\
&   +2 \sin {t} \biggl(\sqrt{2 \pi } S\biggl(\sqrt{\frac{2}{\pi}} \sqrt{t}\biggr) (t+\sin {2t})+2 \sqrt{t} \cos {t}\biggr)\biggr\}.
\end{align*}
\end{lemma}
\proofbegin
The relevant integrals can be computed using the type of argument employed in the proof of lemmaÊ \ref{T1lemma} and integration by parts.
\proofend

\begin{lemma}
The function $I_{shh}(t)$ defined in (\ref{Ishhdef}) can be expressed in the following form:
\begin{align} \nonumber
I_{shh} = \;& 2A(H_s h) + 2\cos(t) I[h A (H_c \cos)] + 2 \sin(t) I[h A(H_s \cos)]
	\\\nonumber
& + 2 \cos(t)I[h A( H_s \sin)]- 2 \sin(t)I[hA (H_c \sin)]	
	\\ \nonumber
&-2A[\cos(t) I(h H_c \cos)]
-2A[\sin(t) I(h H_s \cos)]
	\\ \nonumber
& -2A[\cos(t) I(h H_s \sin)]
+2A[\sin(t) I(h H_c \sin)]
	\\ \nonumber
& -2H_c A(H_c \cos) -2H_cA(H_s \sin)
+ 2A(H_c^2 \cos) 
	\\ \label{Ishhexpression}
& -2H_s A(H_s \cos) + 2H_sA(H_c \sin)
+ 2A(H_s^2 \cos).
\end{align}
\end{lemma}
\proofbegin
Exchanging the order of integration in the definition (\ref{Ishhdef}) of $I_{shh}$ and then computing the integral with respect to $dt'$, we find
\begin{align} \nonumber
I_{shh} =\;& \int_0^t dt'' h(t'') \int_0^{t''}dt''' h(t''')\int_{t''}^t dt' \frac{\sin{t'}}{(t - t' + t'' - t''')^{3/2}}
  	\\ \label{Ishh1}
= & \int_0^t dt'' h(t'') \int_0^{t''} dt''' h(t''')\biggl\{ \frac{2\sin{t}}{\sqrt{t'' - t'''}} - \frac{2\sin{t''}}{\sqrt{t - t'''}}
	\\ \nonumber
& + 2\sqrt{2\pi}\cos(t+ t'' - t''')\biggl[C\biggl(\sqrt{\frac{2}{\pi}}\sqrt{t'' - t'''}\biggr) - C\biggl(\sqrt{\frac{2}{\pi}}\sqrt{t - t'''}\biggr)\biggr]
	\\ \nonumber
& + 2\sqrt{2\pi}\sin(t+ t'' - t''')\biggl[S\biggl(\sqrt{\frac{2}{\pi}}\sqrt{t'' - t'''}\biggr) - S\biggl(\sqrt{\frac{2}{\pi}}\sqrt{t - t'''}\biggr)\biggr]\biggr\}.	
\end{align}
Integration by parts implies that the terms involving $C\bigl(\sqrt{\frac{2}{\pi}}\sqrt{t'' - t'''}\bigr)$ and $S\bigl(\sqrt{\frac{2}{\pi}}\sqrt{t'' - t'''}\bigr)$ can be written as
$$2\int_0^t dt'' h(t'') \int_0^{t''}dt'''  \biggl(\int_0^{t'''} h(r) \cos(t+ t'' - r)dr\biggr) \frac{\cos(t'' - t''')}{\sqrt{t'' - t'''}}$$
and
$$2\int_0^t dt'' h(t'') \int_0^{t''}dt'''  \biggl(\int_0^{t'''} h(r) \sin(t+ t'' - r)dr\biggr) \frac{\sin(t'' - t''')}{\sqrt{t'' - t'''}},$$
respectively. Using the identity
$$\cos(t + t'' - r) \cos(t'' - t''') + \sin(t + t'' - r) \sin(t'' - t''')
= \cos(r - t - t'''),$$
it follows that the total contribution of these two terms is given by
$$2 \int_0^t dt'' h(t'') \int_0^{t''}dt'''  \biggl(\int_0^{t'''} h(r) \cos(r - t - t''')dr\biggr) \frac{1}{\sqrt{t'' - t'''}}.$$

Similarly, integrating by parts twice, we infer that the terms on the rhs of (\ref{Ishh1}) involving $C\bigl(\sqrt{\frac{2}{\pi}}\sqrt{t - t'''}\bigr)$ and $S\bigl(\sqrt{\frac{2}{\pi}}\sqrt{t - t'''}\bigr)$ can be written as
\begin{align*}
& -2\int_0^tdt'' \biggl(\int_0^{t''} h(\tau) \biggl(\int_0^{\tau} h(r)\cos(t + \tau - r) dr\biggr)d\tau\biggr) \frac{\cos(t- t'')}{\sqrt{t - t''}}
	\\
&  -2 \int_0^tdt'' h(t'') \int_0^{t''} dt''' \biggl(\int_0^{t'''} h(r)\cos(t + t'' - r) dr\biggr)\frac{\cos(t - t''')}{\sqrt{t - t'''}}
\end{align*}
and
\begin{align*}
& -2\int_0^tdt'' \biggl(\int_0^{t''} h(\tau) \biggl(\int_0^{\tau} h(r)\sin(t + \tau - r) dr\biggr)d\tau\biggr) \frac{\sin(t- t'')}{\sqrt{t - t''}}
	\\
&  -2 \int_0^tdt'' h(t'') \int_0^{t''} dt''' \biggl(\int_0^{t'''} h(r)\sin(t + t'' - r) dr\biggr)\frac{\sin(t - t''')}{\sqrt{t - t'''}},
\end{align*}
respectively. Using the identities
\begin{align*}
& \cos(t + \tau - r) \cos(t- t'') + \sin(t + \tau - r) \sin(t- t'')
= \cos(r - t'' - \tau),
	\\
&\cos(t + t'' - r) \cos(t - t''') + \sin(t + t'' - r) \sin(t - t''')
= \cos(r - t'' - t'''),
\end{align*}
we infer that the total contribution of these two terms is given by
\begin{align*}
&-2\int_0^tdt'' \biggl(\int_0^{t''} h(\tau) \biggl(\int_0^{\tau} h(r)\cos(r - t'' - \tau) dr\biggr)d\tau\biggr) \frac{1}{\sqrt{t - t''}}
	\\
&  -2 \int_0^tdt'' h(t'') \int_0^{t''} dt''' \biggl(\int_0^{t'''} h(r)\cos(r - t'' - t''') dr\biggr)\frac{1}{\sqrt{t - t'''}}.
\end{align*}
In summary,
\begin{align} \label{Ishhsummary}
I_{shh} = & \int_0^t dt'' h(t'') \int_0^{t''} dt''' h(t''')\biggl( \frac{2\sin{t}}{\sqrt{t'' - t'''}} - \frac{2\sin{t''}}{\sqrt{t - t'''}}\biggr)
	\\ \nonumber
&+ 2 \int_0^t dt'' h(t'') \int_0^{t''}dt'''  \left(\int_0^{t'''} h(r) \cos(r - t - t''')dr\right) \frac{1}{\sqrt{t'' - t'''}}
	\\ \nonumber
&-2\int_0^tdt'' \biggl(\int_0^{t''} h(\tau) \biggl(\int_0^{\tau} h(r)\cos(r - t'' - \tau) dr\biggr)d\tau\biggr) \frac{1}{\sqrt{t - t''}}
	\\ \nonumber
&  -2 \int_0^tdt'' h(t'') \int_0^{t''} dt''' \biggl(\int_0^{t'''} h(r)\cos(r - t'' - t''') dr\biggr)\frac{1}{\sqrt{t - t'''}}.
\end{align}
The term involving $\frac{2\sin{t}}{\sqrt{t'' - t'''}}$ equals 
$$2\sin(t) I(hAh) = 2\pi\sin(t)H_s.$$ 
Integrating by parts with respect to $dt''$, we can write the term involving $-\frac{2\sin{t''}}{\sqrt{t - t'''}}$ as
\begin{align*}
-2H_s(t)(Ah)(t) + 2(A(H_s h))(t) = -2\pi H_s(t) \sin(t) + 2(A(H_s h))(t).
\end{align*}
Moreover, using that 
$$\cos(r-t'' - t''') = \cos{t''}\cos(r-t''') + \sin{t''}\sin(r-t''')$$ 
and integrating by parts with respect to $dt''$, we can write the last line of (\ref{Ishhsummary}) in the following form:
\begin{align*}
& -2 H_c(t) \int_0^{t} dt''' \biggl(\int_0^{t'''} h(r)\cos(r - t''') dr\biggr)\frac{1}{\sqrt{t - t'''}}
  	\\
& +2 \int_0^t dt'' H_c(t'') \int_0^{t''} h(r)\cos(r - t'') dr\frac{1}{\sqrt{t - t''}}
	\\
& -2 H_s(t) \int_0^{t} dt''' \biggl(\int_0^{t'''} h(r)\sin(r - t''') dr\biggr)\frac{1}{\sqrt{t - t'''}}
  	\\
& +2 \int_0^t dt'' H_s(t'') \int_0^{t''} h(r)\sin(r - t'') dr\frac{1}{\sqrt{t - t''}}
	\\
= & -2H_c(t)A(H_c \cos)(t) -2H_c(t)A(H_s \sin)(t)
+ 2A(H_c^2 \cos)(t) + 2A(H_c H_s \sin)(t)
	\\
& -2H_s(t)A(H_s \cos)(t) + 2H_s(t)A(H_c \sin)(t)
+ 2A(H_s^2 \cos)(t) - 2A(H_s H_c \sin)(t).
\end{align*}
Expanding the cosine functions in the second and third lines of (\ref{Ishhsummary}) in a similar way, we find (\ref{Ishhexpression}).
\proofend

\begin{lemma}
The function $I_{hsh}(t)$ defined in (\ref{Ihshdef}) can be expressed in the following form:
\begin{align}\nonumber
I_{hsh} = &\;4A(hH_s) - 2\pi H_s\sin{t}  - 3 A(H_c^2\cos) - 6A(H_cH_s\sin) + 3A(H_s^2 \cos)
	\\\nonumber
&+ 2A[\cos(t)I(hH_c \cos )] + 2A[\sin(t) I(hH_s\cos)] 
	\\\label{Ihshexpression}
& - 2A[\cos(t)I(hH_s \sin)] + 2A[\sin(t)I(hH_c\sin)]  
	\\\nonumber
&+2H_cA(H_c\cos) + 2H_cA(H_s\sin) + 2H_sA(H_c\sin) - 2H_sA(H_s\cos).
\end{align}
\end{lemma}
\proofbegin
Exchanging the order of integration in the definition (\ref{Ihshdef}) of $I_{hsh}$ and then computing the integral with respect to $dt''$, we find
\begin{align*}
I_{hsh}  = & \int_0^t dt' h(t')  \int_0^{t'} dt''' h(t''') \int_{t'''}^{t'} dt'' \frac{\sin{t''}}{(t - t' + t'' - t''')^{3/2}}
	\\
= & 	\int_0^t dt' h(t')  \int_0^{t'} dt''' h(t''')\biggl(\frac{2\sin{t'''}}{\sqrt{t - t'}} - \frac{2\sin{t'}}{\sqrt{t - t'''}}\biggr)
	\\
&+ 	2 \sqrt{2 \pi } \int_0^t dt' h(t')  \int_0^{t'} dt''' h(t''')\cos (t-t'-t''') 
	\\
& \hspace{4cm}\times \biggl(C\biggl(\sqrt{\frac{2}{\pi }}
   \sqrt{t-t'''}\biggr)-C\biggl(\sqrt{\frac{2}{\pi }} \sqrt{t-t'}\biggr)\biggr) 
   	\\
& +2 \sqrt{2 \pi } \int_0^t dt' h(t')  \int_0^{t'} dt''' h(t''')\sin (t-t'-t''')
	\\
& \hspace{4cm}\times \biggl(S\biggl(\sqrt{\frac{2}{\pi }} \sqrt{t-t'''}\biggr)-S\biggl(\sqrt{\frac{2}{\pi }}\sqrt{t-t'}\biggr)\biggr). 
\end{align*}
Exchanging the order of integration, integrating by parts with respect to $dt'$, and then changing back the order of integration, we find that the terms involving $C\bigl(\sqrt{\frac{2}{\pi}}\sqrt{t - t'}\bigr)$ and $S\bigl(\sqrt{\frac{2}{\pi}}\sqrt{t - t'}\bigr)$ can be written as
$$-2\int_0^t dt' \int_0^{t'}dt''' h(t''')\left(\int_{t'''}^{t'} h(r) \cos(t - r - t''')dr\right) \frac{\cos(t - t')}{\sqrt{t - t'}}$$
and
$$-2\int_0^t dt' \int_0^{t'}dt''' h(t''')\left(\int_{t'''}^{t'} h(r) \sin(t - r - t''')dr\right) \frac{\sin(t - t')}{\sqrt{t - t'}},$$
respectively. Using the identity
$$\cos(t - r - t''')\cos(t - t') +  \sin(t - r - t''')\sin(t - t')
= \cos(r - t' + t'''),$$
it follows that the total contribution of these two terms is given by
$$-2\int_0^t \frac{dt'}{\sqrt{t - t'}} \int_0^{t'}dt''' h(t''')\int_{t'''}^{t'} dr h(r) \cos(r - t' + t''').$$

Similarly, integrating by parts twice, we infer that the terms involving $C(\sqrt{\frac{2}{\pi}}\sqrt{t - t'''})$ and $S(\sqrt{\frac{2}{\pi}}\sqrt{t - t'''})$ can be written as
\begin{align*}
&2\int_0^tdt' \biggl(\int_0^{t'} h(\tau) \biggl(\int_0^{\tau} h(r)\cos(t - \tau - r) dr\biggr)d\tau\biggr) \frac{\cos(t- t')}{\sqrt{t - t'}}
	\\
&+2 \int_0^tdt' h(t') \int_0^{t'} dt''' \biggl(\int_0^{t'''} h(r)\cos(t - t' - r) dr\biggr)\frac{\cos(t - t''')}{\sqrt{t - t'''}}
\end{align*}
and
\begin{align*}
&2\int_0^tdt' \biggl(\int_0^{t'} h(\tau) \biggl(\int_0^{\tau} h(r)\sin(t - \tau - r) dr\biggr)d\tau\biggr) \frac{\sin(t- t')}{\sqrt{t - t'}}
	\\
&+2 \int_0^tdt' h(t') \int_0^{t'} dt''' \biggl(\int_0^{t'''} h(r)\sin(t - t' - r) dr\biggr)\frac{\sin(t - t''')}{\sqrt{t - t'''}},
\end{align*}
respectively. Using the identities
\begin{align*}
&\cos(t - \tau - r) \cos(t - t') + \sin(t - \tau - r) \sin(t- t')
= \cos(r - t' + \tau),
	\\
&\cos(t - t' - r)\cos(t - t''') + \sin(t - t' - r)\sin(t - t''')
= \cos(t' - t''' + r),
\end{align*}
we infer that the total contribution of these two terms is
\begin{align*}
&2\int_0^tdt' \biggl(\int_0^{t'} h(\tau) \biggl(\int_0^{\tau} h(r) \cos(r - t' + \tau) dr\biggr)d\tau\biggr) \frac{1}{\sqrt{t - t'}}
	\\
&+2 \int_0^tdt' h(t') \int_0^{t'} dt''' \biggl(\int_0^{t'''} h(r)\cos(t' - t''' + r)dr\biggr)\frac{1}{\sqrt{t - t'''}}.
\end{align*}
In summary,
\begin{align}\label{Ihshsummary}
I_{hsh} &= 	\int_0^t dt' h(t')  \int_0^{t'} dt''' h(t''')\biggl(\frac{2\sin{t'''}}{\sqrt{t - t'}} - \frac{2\sin{t'}}{\sqrt{t - t'''}}\biggr)
	\\\nonumber
& -2\int_0^t dt' \int_0^{t'}dt''' h(t''')\left(\int_{t'''}^{t'} h(r) \cos(r - t' + t''')dr\right) \frac{1}{\sqrt{t - t'}}
	\\\nonumber
&+ 2\int_0^tdt' \biggl(\int_0^{t'} h(\tau) \biggl(\int_0^{\tau} h(r) \cos(r - t' + \tau) dr\biggr)d\tau\biggr) \frac{1}{\sqrt{t - t'}}
	\\\nonumber
&+2 \int_0^tdt' h(t') \int_0^{t'} dt''' \biggl(\int_0^{t'''} h(r)\cos(t' - t''' + r)dr\biggr)\frac{1}{\sqrt{t - t'''}}.
\end{align}
Integration by parts with respect to $dt'$ implies that the term involving $ \frac{2\sin{t'}}{\sqrt{t - t'''}}$ equals
$$-2H_s(t)Ah(t) + 2A(hH_s)(t) = -2\pi H_s(t)\sin{t} + 2A(hH_s)(t).$$
We next consider the second line of (\ref{Ihshsummary}): the integrand in the factor
$$\int_0^{t'}dt''' h(t''')\int_{t'''}^{t'} dr  h(r) \cos(r - t' + t''')$$
is symmetric in $r$ and $t'''$, therefore we can write this factor as
$$\frac{1}{2} \int_0^{t'}dt''' \int_{0}^{t'} drh(t''') h(r) \cos(r - t' + t''').$$
Expanding $\cos(r - t' + t''')$, we see that the second line can be written as
$$-A(H_c^2\cos) - 2A(H_cH_s\sin) + A(H_s^2 \cos).$$
The third line can be handled by expanding $\cos(r - t' + \tau)$.
Finally, using the identity 
$$\cos(t' - t''' + r) = \cos{t'}\cos(t''' - r) + \sin{t'}\sin(t''' - r)$$ 
and integrating by parts with respect to $dt'$, we can write the last line of (\ref{Ihshsummary}) in the following form:
\begin{align*}
 & 2H_c(t)\int_0^t \frac{dt'''}{\sqrt{t - t'''}}\int_0^{t'''}drh(r)\cos(t''' - r)
  - 2\int_0^t dt' H_c(t') \int_0^{t'}dr h(r) \frac{\cos(t' - r)}{\sqrt{t - t'}}
  	\\
&  +  2H_s(t)\int_0^t \frac{dt'''}{\sqrt{t - t'''}}\int_0^{t'''}drh(r)\sin(t''' - r)
  - 2\int_0^t dt' H_s(t') \int_0^{t'}dr h(r) \frac{\sin(t' - r)}{\sqrt{t - t'}}
  	\\
&= 2H_cA(H_c\cos) + 2H_cA(H_s\sin) - 2A(H_c^2\cos) - 2A(H_cH_s\sin)
	\\
& + 2H_sA(H_c\sin) - 2H_sA(H_s\cos) - 2A(H_cH_s\sin ) + 2A(H_s^2 \cos ).
\end{align*}
Summing up the various contributions, we find (\ref{Ihshexpression}).
\proofend

\begin{lemma}
The function $I_{hhs}(t)$ defined in (\ref{Ihhsdef}) can be expressed in the following form:
\begin{align} \nonumber
I_{hhs} = &\; 2A(hH_s) + 2A(H_c^2\cos) - 2A[\cos(t) I(h H_c \cos)] - 2A[\sin(t)I(h H_s \cos)] 
	\\\nonumber
&+ 2A(H_s^2 \cos) - 2A[\cos(t) I(h H_s \sin)] + 2A[\sin(t) I(h H_c\sin)]
	\\\nonumber
&+2\sqrt{2\pi}C\biggl(\sqrt{\frac{2}{\pi}}\sqrt{t}\biggr)\Bigl[-\cos(t) I(h H_c \cos)
+ \sin(t) I(h H_s \cos)
	\\\nonumber
&\hspace{3.9cm}-\cos(t) I(hH_s \sin) -\sin(t) I(h H_c\sin)\Bigr]
	\\\nonumber
&+2\sqrt{2\pi}S\biggl(\sqrt{\frac{2}{\pi}}\sqrt{t}\biggr)\Bigl[-\cos(t) I(h H_s \cos)
- \sin(t) I(hH_c \cos)
	\\\nonumber
&\hspace{3.9cm}+ \cos(t) I(h H_c\sin) - \sin(t) I(hH_s \sin)\Bigr]
	\\ \label{Ihhsexpression}
& + 2\int_0^t dt' h(t') \int_0^{t'}dt''  \frac{\cos{t''}H_c(t'') + \sin{t''}H_s(t'')}{\sqrt{ t- t' + t''}}	.
\end{align}
\end{lemma}
\proofbegin
Computing the integral with respect to $dt'''$ in the definition (\ref{Ihhsdef}) of $I_{hhs}$, we find
\begin{align*}
I_{hhs} = &\int_0^t dt' h(t')  \int_0^{t'} dt'' h(t'')\frac{2\sin{t''}}{\sqrt{t - t'}}
	\\
&+ 2 \sqrt{2 \pi } \int_0^t dt' h(t')  \int_0^{t'} dt'' h(t'')\cos (t-t'+t'') 
	\\
&\hspace{3cm} \times \biggl(C\biggl(\sqrt{\frac{2}{\pi }}
   \sqrt{t-t'}\biggr)-C\biggl(\sqrt{\frac{2}{\pi }} \sqrt{t-t' + t''}\biggr)\biggr) 
   	\\
& +2 \sqrt{2 \pi } \int_0^t dt' h(t')  \int_0^{t'} dt'' h(t'')\sin(t-t'+t'')
	\\
&\hspace{3cm} \times\biggl(S\biggl(\sqrt{\frac{2}{\pi }} \sqrt{t-t'}\biggr)-S\biggl(\sqrt{\frac{2}{\pi }}\sqrt{t-t'+t''}\biggr)\biggr).
\end{align*}
Exchanging the order of integration, integrating by parts with respect to $dt'$, and then changing back the order of integration, we find that the terms involving $C\bigl(\sqrt{\frac{2}{\pi}}\sqrt{t - t'}\bigr)$ and $S\bigl(\sqrt{\frac{2}{\pi}}\sqrt{t - t'}\bigr)$, can be written as
$$2\int_0^t dt' \int_0^{t'}dt'' h(t'') \biggl(\int_{t''}^{t'} h(r)\cos(t - r + t'')dr\biggr) \frac{\cos(t - t')}{\sqrt{t - t'}} $$
and
$$2\int_0^t dt' \int_0^{t'}dt'' h(t'') \biggl(\int_{t''}^{t'} h(r)\sin(t - r + t'')dr\biggr) \frac{\sin(t - t')}{\sqrt{t - t'}},$$
respectively. Using the identity
\begin{align}\nonumber
\cos(t - r + t'')\cos(t - t') + \sin(t - r + t'')\sin(t - t') = \cos(r - t' - t''),
\end{align}
we see that the total contribution of these two terms is
$$2\int_0^t dt' \int_0^{t'}dt'' h(t'') \biggl(\int_{t''}^{t'} h(r)\cos(r - t' - t'')dr\biggr) \frac{1}{\sqrt{t - t'}}.$$

Similarly, integrating by parts with respect to $dt''$, the terms involving $C\bigl(\sqrt{\frac{2}{\pi}}\sqrt{t - t' + t''}\bigr)$ and $S\bigl(\sqrt{\frac{2}{\pi}}\sqrt{t - t' + t''}\bigr)$ can be written as
\begin{align*}
& -2\sqrt{2\pi}C\biggl(\sqrt{\frac{2}{\pi}}\sqrt{t}\biggr) \int_0^t dt' h(t') \int_0^{t'} h(r) \cos(t - t' + r) dr
	\\
& + 2\int_0^t dt' h(t') \int_0^{t'}dt'' \biggl(\int_0^{t''} h(r) \cos(t - t' + r) dr\biggr) \frac{\cos(t - t' + t'')}{\sqrt{ t- t' + t''}}	
\end{align*}
and
\begin{align*}
& -2\sqrt{2\pi}S\biggl(\sqrt{\frac{2}{\pi}}\sqrt{t}\biggr) \int_0^t dt' h(t') \int_0^{t'} h(r) \sin(t - t' + r) dr
	\\
& + 2\int_0^t dt' h(t') \int_0^{t'}dt'' \biggl(\int_0^{t''} h(r) \sin(t - t' + r) dr\biggr) \frac{\sin(t - t' + t'')}{\sqrt{ t- t' + t''}}.	
\end{align*}
Using the identity
\begin{align}\nonumber
\cos(t - t' + r)\cos(t - t' + t'') + \sin(t - t' + r)\sin(t - t' + t'') = \cos(r - t''),
\end{align}
we infer that the total contribution of these two terms is
\begin{align*}
& -2\sqrt{2\pi}C\biggl(\sqrt{\frac{2}{\pi}}\sqrt{t}\biggr) \int_0^t dt' h(t') \int_0^{t'} h(r) \cos(t - t' + r) dr
	\\
& -2\sqrt{2\pi}S\biggl(\sqrt{\frac{2}{\pi}}\sqrt{t}\biggr) \int_0^t dt' h(t') \int_0^{t'} h(r) \sin(t - t' + r) dr	
	\\
& + 2\int_0^t dt' h(t') \int_0^{t'}dt'' \biggl(\int_0^{t''} h(r) \cos(r - t'') dr\biggr) \frac{1}{\sqrt{ t- t' + t''}}.	
\end{align*}
In summary, 
\begin{align*}
I_{hhs} & = \int_0^t dt' h(t')  \int_0^{t'} dt'' h(t'')\frac{2\sin{t''}}{\sqrt{t - t'}}
	\\
& + 2\int_0^t dt' \int_0^{t'}dt'' h(t'') \biggl(\int_{t''}^{t'} h(r)\cos(r - t' - t'')dr\biggr) \frac{1}{\sqrt{t - t'}}
	\\
& -2\sqrt{2\pi}C\biggl(\sqrt{\frac{2}{\pi}}\sqrt{t}\biggr) \int_0^t dt' h(t') \int_0^{t'} h(r) \cos(t - t' + r) dr
	\\
& -2\sqrt{2\pi}S\biggl(\sqrt{\frac{2}{\pi}}\sqrt{t}\biggr) \int_0^t dt' h(t') \int_0^{t'} h(r) \sin(t - t' + r) dr	
	\\
& + 2\int_0^t dt' h(t') \int_0^{t'}dt'' \biggl(\int_0^{t''} h(r) \cos(r - t'') dr\biggr) \frac{1}{\sqrt{ t- t' + t''}}	.
\end{align*}
The first term on the rhs equals $2A(hH_s)$. Expanding the trigonometric functions in the remaining terms, we arrive at (\ref{Ihhsexpression}).
\proofend

We will find the asymptotics of $g_{13}$ by determining the asymptotics of the various terms appearing in the above lemmas. We begin with the functions $H_c$ and $H_s$, which can be expressed in terms of $C$ and $S$.

\begin{lemma}\label{HcHslemma}
The functions $H_c$ and $H_s$ defined in (\ref{HcHsdef}) satisfy the following estimates:
\begin{align*}
 H_c(t) & =  \frac{1}{4} \biggl\{\sqrt{2 \pi } C\biggl(\sqrt{\frac{2}{\pi }} \sqrt{t}\biggr) (2 t+\sin {2t})-\sqrt{2 \pi }
   S\biggl(\sqrt{\frac{2}{\pi }} \sqrt{t}\biggr) \cos {2t}-2 \sqrt{t} \sin (t)\biggr\}
	\\
& = \frac{1}{2} \sqrt{\frac{\pi }{2}} t+\frac{1}{4} \sqrt{\frac{\pi }{2}} (\sin {2t}-\cos {2t})
 -\frac{\sin (t)}{2 t^{3/2}}  +\frac{3\cos (t)}{4 t^{5/2}}+\frac{15 \sin (t)}{4 t^{7/2}} + O(t^{-9/2}).
	\\
 H_s(t) & = \frac{1}{4} \biggl\{2 \sqrt{t} \cos(t)-\sqrt{2 \pi }\biggl(C\biggl(\sqrt{\frac{2}{\pi }} \sqrt{t}\biggr) \cos (2t)+S\biggl(\sqrt{\frac{2}{\pi }} \sqrt{t}\biggr) (\sin {2t}-2 t)\biggr)\biggr\}
   	\\
&=\frac{1}{2} \sqrt{\frac{\pi }{2}} t-\frac{1}{4}
   \sqrt{\frac{\pi }{2}} (\sin {2t}+\cos {2t})
+\frac{\cos (t)}{2 t^{3/2}} +\frac{3\sin (t)}{4 t^{5/2}}-\frac{15 \cos (t)}{4 t^{7/2}} + O(t^{-9/2}).
\end{align*}
\end{lemma}
\proofbegin
This is a direct consequence of using the expression (\ref{hintermsofCS}) for $h(t)$. 
\proofend
 
Since $H_c$ and $H_s$ are of $O(t)$ as $tÊ\to \infty$, by scrutinizing the above lemmas it becomes clear that in order to determine the asymptotics of $g_{13}$ to order $O(t^{-1/16})$, we need to find the asymptotics of the following terms to $O(t^{-1/16})$:
\begin{align}\nonumber
& A(H_s h), \quad 
A[\cos(t) I(h H_c \cos)], \quad A[\sin(t) I(h H_s \cos)], \quad A[\cos(t) I(h H_s \sin)], 
	\\  \nonumber
&A[\sin(t) I(h H_c \sin)], \quad A(H_c^2 \cos), \quad A(H_s^2 \cos), \quad A(H_cH_s \sin),
	\\\nonumber
& I[h A (H_c \cos)], \quad I[h A(H_s \cos)], \quad I[h A( H_s \sin)], \quad I[hA (H_c \sin)],
	\\ \label{listoforder18}
& I(h H_c\cos), \quad  I(h H_s\cos), \quad  I(h H_s\sin), \quad
  I(h H_c\sin),
\end{align}
as well as of the following terms to $O(t^{-17/16})$:
\begin{align}\label{listoforder98}
  A(H_c\cos), \quad A(H_s\sin), \quad A(H_c\sin), \quad A(H_s\cos).
\end{align}
The following lemma will be the main tool for finding the asymptotics of the relevant Abel integrals.

\begin{lemma}\label{abelasymptoticslemma}
Let $N >0$. Let $f(t)$, $t \geq 0$, be a smooth function which is bounded as $t\to 0$, and which satisfies
\begin{align*}
f(t) = f_a(t) + O(t^{-1}),\qquad t \to \infty,
\end{align*}
where
\begin{align*}
f_a(t) = \sum_{j=0}^4 t^{\frac{j}{2}}\left(f_j + \sum_{n = 1}^N (c_{jn} \cos{nt} + s_{jn} \sin{nt})\right)
+ \frac{\hat{f}_1 + \sum_{n = 1}^N (\hat{c}_{1n}\cos{nt} + \hat{s}_{1n}\sin{nt})}{\sqrt{t}} ,
\end{align*}
and $f_j$, $c_{jn}$, $s_{jn}$, $\hat{c}_{1n}$, $\hat{s}_{1n}$, $j=0,\dots,4$, $n = 1, \dots,N$, and $\hat{f}_1$ are constants. Then, the asymptotic behavior of $Af(t)$ as $t\to \infty$ is given by
\begin{align*}
  Af(t) = \;&
  \frac{\sqrt{\pi } f_j t^{\frac{j + 1}{2}} \Gamma \left(\frac{j}{2}+1\right)}{\Gamma \left(\frac{j+3}{2}\right)} + \hat{f}_1\pi
 	\\  
& +  \sqrt{\frac{\pi }{2}}   \sum_{j=0}^4 \sum_{n = 1}^N \frac{c_{jn}}{32 n^{5/2}}   t^{\frac{j}{2}-2} 
\biggl\{\left(-3 (j-2) j -8 j n t +32 n^2 t^2\right) \sin {nt}
 	\\
&\hspace{5cm} +\left(-3 (j-2) j +8 j n t+32 n^2 t^2\right) \cos {nt}\biggr\} 
	\\
& + \sqrt{\frac{\pi }{2}}  \sum_{j=0}^4 \sum_{n = 1}^N  \frac{s_{jn}}{4 n^{3/2}} t^{\frac{j}{2}-1} ((j+4 n t) \sin
   {nt}+(j-4 n t) \cos {nt}) + O(t^{-1/16}).
\end{align*}
\end{lemma}
\proofbegin
Since $f(t) = O(t^2)$ as $t\to \infty$, there exist constants $K>0$ and $M>0$ such that
\begin{align}\label{int0t18}
\left|\int_0^{t^{1/8}} \frac{f(s)}{\sqrt{t-s}} ds\right|
\leq \int_0^{K} \frac{|f(s)|}{\sqrt{t-s}} ds + M \int_K^{t^{1/8}} \frac{s^2}{\sqrt{t-s}} ds, \qquad t^{1/8} > K.
\end{align}
The first integral on the rhs is of $O(t^{-1/2})$ as $t \to \infty$, while the second can be computed explicitly and is of $O(t^{-1/8})$ as $t \to \infty$. 

On the other hand, by assumption, there exists $L> 0$ such that $|f(t) - f_a(t)| < \frac{L}{t}$ for $t$ large enough.
Thus, as $t\to \infty$,
\begin{align*}
& \left|\int_{t^{1/8}}^t \frac{f(s) - f_a(s)}{\sqrt{t-s}} ds\right|
\leq L \int_{t^{1/8}}^t \frac{1}{s\sqrt{t-s}} ds = O(t^{-1/4}).
\end{align*}
We conclude that to $O(t^{-1/8})$ the asymptotics of $Af$ coincide with the asymptotics of $\int_{t^{1/8}}^t \frac{f_a(s)}{\sqrt{t - s}} ds$.
In view of the estimates
\begin{align*}
& \left|\int_{t^{1/8}}^t \frac{\sin{ns}}{\sqrt{s}\sqrt{t-s}} ds\right|
\leq  \frac{1}{t^{1/16}} \int_{t^{1/8}}^t \frac{\sin{ns}}{\sqrt{t-s}} ds
= O(t^{-1/16}),
	\\
& \left|\int_{t^{1/8}}^t \frac{\cos{ns}}{\sqrt{s}\sqrt{t-s}} ds\right|
\leq  \frac{1}{t^{1/16}} \int_{t^{1/8}}^t \frac{\cos{ns}}{\sqrt{t-s}} ds
= O(t^{-1/16}),
\end{align*}
the contribution to the asymptotics of the terms involving $\hat{c}_{1n}$ and $\hat{s}_{1n}$ are of $O(t^{-1/16})$. 
The contribution of the term involving $\hat{f}_1$ is $\hat{f}_1\pi + O(t^{-1/16})$. The contributions of the remaining terms can be obtained by expressing the integrals\footnote{The lower limit of integration here can be taken to be $t^{1/8}$ or $0$; an estimate as in (\ref{int0t18}) shows that the difference is of $O(t^{-1/8})$.} 
$$\int_0^t t^{\frac{j}{2}}\left(f_j + \sum_{n = 1}^N (c_{jn} \cos{nt} + s_{jn} \sin{nt})\right), \qquad j = 0, \dots, 4,$$
in terms of hypergeometric functions and computing the asymptotics.
\proofend

Although lemma \ref{abelasymptoticslemma} is sufficient for most purposes, we need a more sophisticated argument to deal with the higher order asymptotics of the terms in (\ref{listoforder98}).

\begin{lemma}\label{98lemma}
Let $A, H_c, H_s$ be defined in (\ref{Adef}) and (\ref{HcHsdef}). Then,
\begin{align*}
&  A(H_c \cos)(t) =\frac{\pi t}{4}(\sin{t} + \cos{t}) -\frac{\pi \sin{t}}{8} -\frac{\pi  \cos{3 t}}{8 \sqrt{3}} - \frac{\sqrt{\pi }}{6 \sqrt{t}}+ O(t^{-9/8}),
  	\\
& A(H_s\sin)(t) =  \frac{\pi  t}{4} (\sin{t} - \cos{t}) +\frac{\pi \sin{t}}{8} +\frac{\pi  \cos{3 t}}{8\sqrt{3}} + \frac{\sqrt{\pi}}{6 \sqrt{t}}+ O(t^{-9/8}),
 	\\
& A(H_c\sin)(t) = \frac{\pi  t }{4}(\sin{t} - \cos{t}) + \frac{\pi  \sin{t}}{4} +\frac{\pi  \cos{t}}{8} -\frac{\pi  \sin {3 t}}{8 \sqrt{3}}
  -\frac{\sqrt{\pi }}{6 \sqrt{t}}  + O(t^{-9/8}),
 	\\
& A(H_s\cos)(t) = \frac{\pi  t }{4}(\sin{t} + \cos{t}) - \frac{\pi  \sin{t}}{4} +\frac{\pi  \cos{t}}{8} -\frac{\pi  \sin {3 t}}{8 \sqrt{3}}
  -\frac{\sqrt{\pi }}{6 \sqrt{t}}  + O(t^{-9/8}).
\end{align*}
\end{lemma}
\proofbegin
Suppose that $f(t)$, $t \geq 0$, is a smooth function which is bounded near $t =0$. 
Let $F(t) = \int_0^t f(s) ds$ and suppose that $f_a(t)$ and $F_a(t)$ are smooth functions bounded near $t = 0$ such that
\begin{align}\label{fFassumptions}
f(t) = f_a(t) + O(t^{-4}), \qquad F(t) = F_a(t) + O(1), \qquad t \to \infty.
\end{align}
Integrating by parts the first term on the rhs of the equation
\begin{align*}
  Af(t) = \int_0^{t^{1/4}} \frac{f(s)ds}{\sqrt{t - s}} + \int_{t^{1/4}}^t \frac{f(s)ds}{\sqrt{t - s}},
\end{align*}
we find
\begin{align*}
  Af(t) = \frac{F(t^{1/4})}{\sqrt{t - t^{1/4}}} - \int_0^{t^{1/4}} \frac{F(s)ds}{2(t - s)^{3/2}} + \int_{t^{1/4}}^t \frac{(f(s) - f_a(s))ds}{\sqrt{t - s}} + \int_{t^{1/4}}^t \frac{f_a(s)ds}{\sqrt{t - s}}.
\end{align*}
Using the estimates
$$\int_0^{t^{1/4}} \frac{ds}{(t-s)^{3/2}} = O(t^{-5/4}), \qquad
\int_{t^{1/4}}^t \frac{ds}{s^4\sqrt{t-s}} = O(t^{-5/4}), \qquad t \to \infty,$$
we conclude that
\begin{align}\label{Afasymptotics}
  Af(t) = \frac{F(t^{1/4})}{\sqrt{t - t^{1/4}}} - \int_0^{t^{1/4}} \frac{F_a(s)ds}{2(t - s)^{3/2}} + \int_{t^{1/4}}^t \frac{f_a(s)ds}{\sqrt{t - s}} + O(t^{-5/4}), \qquad t \to \infty.
\end{align}

Now suppose that $f(t) = H_s(t) \sin{t}$. Then equation (\ref{fFassumptions}) is satisfied with
$$F_a(t) = -\frac{1}{2} \sqrt{\frac{\pi }{2}} t \cos{t}$$
and
\begin{align*}
f_a(t) = &\;\frac{1}{2}\sqrt{\frac{\pi }{2}} t \sin (t)+\frac{1}{8} \sqrt{\frac{\pi }{2}}\bigl(\sin (t)-\sin (3 t)- \cos (t)+ \cos(3 t)\bigr)
	\\
& +\frac{\sin {2t}}{4 t^{3/2}} 
 +\frac{3}{8 t^{5/2}}
-\frac{3 \cos {2t}}{8 t^{5/2}}
-\frac{15 \sin {2t}}{8 t^{7/2}}.
\end{align*}
Using that
\begin{align*}
F(t) =& -\frac{1}{2} \sqrt{\frac{\pi }{2}} t \cos (t)
+ \frac{\sqrt{\pi }}{6} 
-\frac{1}{24} \sqrt{\frac{\pi }{2}} \bigl(9\sin(t)+ \sin (3 t) - 3 \cos (t) + \cos(3 t)\bigr)
	\\
& -\frac{1}{4 t^{3/2}}
- \frac{\cos (2 t)}{8 t^{3/2}}
- \frac{9 \sin (2 t)}{32 t^{5/2}}
+ O(t^{-7/2}), \qquad t \to \infty,
\end{align*}
it is easy to find the asymptotics of all the terms in (\ref{Afasymptotics}) up to $O(t^{-5/4})$ except for the contribution from the following terms in $f_a$:
\begin{align}\label{threedifficultterms}
\frac{\sin {2t}}{4 t^{3/2}}, \quad -\frac{3 \cos {2t}}{8 t^{5/2}}, \quad -\frac{15 \sin {2t}}{8 t^{7/2}}.
\end{align}
In order to handle the first of these terms, we consider an integral of the form
$$\int_{t^{1/4}}^t \frac{\sin (n s)}{s^{3/2}\sqrt{t - s}}ds, \qquad n > 0.$$
We write
\begin{align}\label{intt14tsinns}
\int_{t^{1/4}}^t \frac{\sin (n s)}{s^{3/2}\sqrt{t - s}}ds
= \int_{t^{1/4}}^{t-1} \frac{\sin (n s)}{s^{3/2}\sqrt{t - s}}ds
+\int_{t-1}^t \frac{\sin (n s)}{s^{3/2}\sqrt{t - s}}ds.
\end{align}
Integrating by parts twice, we can write the first term on the rhs as
\begin{align}\nonumber
& -\frac{\cos (ns)}{n s^{3/2}\sqrt{t - s}}\bigg|_{s=t^{1/4}}^{t-1}
+ \frac{\sin(ns)}{n^2} \frac{4 s-3 t}{2 s^{5/2} (t-s)^{3/2}} \bigg|_{s=t^{1/4}}^{t-1}
	\\ \label{IBPthreeterms}
& - \int_{t^{1/4}}^{t-1} \frac{\sin(ns)}{n^2} \frac{3 \left(8 s^2-12 s t+5 t^2\right)}{4 s^{7/2} (t-s)^{5/2}} ds.
\end{align}
Using that $8 s^2-12 s t+5 t^2 > 0$ for $s \in (t^{1/4}, t-1)$ for $t$ large enough and that
$$ \int_{t^{1/4}}^{t-1} \frac{3 \left(8 s^2-12 s t+5 t^2\right)}{4 s^{7/2} (t-s)^{5/2}} ds = O(t^{-9/8}),$$
we deduce that the last integral in (\ref{IBPthreeterms}) is of $O(t^{-9/8})$.
On the other hand, the second term on the rhs of (\ref{intt14tsinns}) can be estimated as follows:
\begin{align*}
 \biggl|\int_{t-1}^t \frac{\sin (n s)}{4 s^{3/2}\sqrt{t - s}}ds\biggr|
 \leq \frac{1}{(t-1)^{3/2}} \int_{t-1}^t \frac{1}{4\sqrt{t - s}}ds = O(t^{-3/2}), \qquad t \to \infty.
\end{align*}
In summary, 
$$\int_{t^{1/4}}^t \frac{\sin (n s)}{s^{3/2}\sqrt{t - s}}ds
= \frac{\cos \left(n t^{1/4}\right)}{n t^{7/8}} + O(t^{-9/8}).$$
This determines the contribution to equation (\ref{Afasymptotics}) from the first term in (\ref{threedifficultterms}). 
An analogous argument shows that the two other terms in (\ref{threedifficultterms}) yield contributions to equation (\ref{Afasymptotics}) of $O(t^{-9/8})$.
Equation (\ref{Afasymptotics}) now gives the asymptotics of the function $A(H_s \sin)$ as stated in the lemma; the other functions can be analyzed in a similar way. 
\proofend

Combining the known asymptotics of $h$, $H_c$, $H_s$ with lemma \ref{98lemma}, we find the following asymptotics of the last eight quantities in (\ref{listoforder18}).

\begin{lemma}\label{lasteightlemma}
Let $h, I, A, H_c, H_s$ be defined in (\ref{g11formula}) and (\ref{Adef})-(\ref{HcHsdef}). Then the following estimates are valid:
\begin{align*}
 I[h A (H_c \cos)] = \frac{\pi^{3/2}}{\sqrt{2}}\biggl\{&\frac{t^2}{8}-\frac{t}{16} -\frac{t \cos{2t}}{8}-\frac{\sin {2t}}{32 \sqrt{3}}
  +\frac{3\sin {2t}}{32}
 	\\
&-\frac{\sin {4t}}{64 \sqrt{3}}- \frac{\cos {2t}}{32 \sqrt{3}}+\frac{\cos {2t}}{32} +\frac{\cos {4t}}{64 \sqrt{3}} \biggr\} + k_1 + O(t^{-1/8}),
	\\ 
 I[h A(H_s \cos)]=
  \frac{\pi^{3/2}}{\sqrt{2}}\biggl\{& \frac{t^2}{8}-\frac{t}{16}-\frac{t \cos {2t}}{8}
  -\frac{\sin {2t}}{32 \sqrt{3}} +\frac{5\sin {2t}}{32}
  	\\
&  +\frac{\sin {4t}}{64 \sqrt{3}}+\frac{\cos {2t}}{32 \sqrt{3}}+\frac{\cos {2t}}{32} +\frac{\cos{4t}}{64 \sqrt{3}} \biggr\} + k_2 + O(t^{-1/8}),
	\\  
 I[h A( H_s \sin)]= \frac{\pi^{3/2}}{\sqrt{2}}\biggl\{&\frac{t}{16}-\frac{t \sin {2t}}{8} +\frac{\sin {2t}}{32 \sqrt{3}} -\frac{\sin {2t}}{32} +\frac{\sin {4t}}{64 \sqrt{3}}
 	\\
&+\frac{\cos {2t}}{32 \sqrt{3}}-\frac{3\cos {2t}}{32} -\frac{\cos {4t}}{64 \sqrt{3}} \biggr\} + k_3 + O(t^{-1/8}),
	\\ 
 I[h A (H_c \sin)] = \frac{\pi^{3/2}}{\sqrt{2}}\biggl\{& \frac{3 t}{16}-\frac{t \sin {2t}}{8}  -\frac{\sin {2t}}{32 \sqrt{3}} -\frac{\sin {2t}}{32} +\frac{\sin {4t}}{64 \sqrt{3}}
 	\\
& +\frac{\cos {2t}}{32 \sqrt{3}}-\frac{5\cos {2t}}{32}+\frac{\cos {4t}}{64 \sqrt{3}}\biggr\} + k_4 + O(t^{-1/8}),
\end{align*}
and
\begin{align*}
 I(h H_c\cos)= & \frac{\pi  t^2}{16}+\frac{\pi  t (\sin {2t} - \cos {2t})}{16}  -\frac{\pi  \sin{4t}}{64} 
-\frac{\sqrt{\frac{\pi }{2}} \sin {t}}{4 \sqrt{t}} + k_5 + O(t^{-3/2}),
	\\ 
 I(h H_s\cos) = &
 \frac{\pi  t^2}{16}-\frac{\pi  t}{16}+\frac{\pi  t (\sin {2t} - \cos {2t})}{16} +\frac{\pi  \cos {2t}}{16} 
 	\\
& 
 +\frac{\pi  \cos {4t}}{64}   -\frac{\sqrt{\frac{\pi }{2}} \sin {t}}{4\sqrt{t}} + k_6 + O(t^{-3/2}),
	\\ 
 I(h H_s\sin)= &\frac{\pi  t^2}{16}-\frac{ \pi  t (\sin {2t} + \cos{2t})}{16}+\frac{\pi  \sin {4t}}{64} 
 +\frac{\sqrt{\frac{\pi }{2}} \cos {t}}{4 \sqrt{t}} + k_7 + O(t^{-3/2}),
	\\ 
  I(h H_c\sin) = &\frac{\pi  t^2}{16}+\frac{\pi  t}{16} - \frac{\pi t (\sin {2t} + \cos{2t})}{16}-\frac{\pi  \cos {2t}}{16} 
  	\\
&  +\frac{ \pi  \cos {4t}}{64}+\frac{\sqrt{\frac{\pi }{2}} \cos{t}}{4 \sqrt{t}} + k_8 + O(t^{-3/2}),
\end{align*}
where $\{k_j\}_1^8$ are real constants.
\end{lemma}
\proofbegin
We consider $I[h A (H_c \cos)]$; the other terms are treated similarly. By (\ref{hasymptotics}) and lemma \ref{98lemma}, we have
$$h A (H_c \cos) = f(t) + f_r(t),$$
where $f(t)$ is a known function and $f_r(t) = O(t^{-9/8})$. 
Choose $K > 0$ and $M >0$ such that
$$|f_r(t)| \leq \frac{M}{t^{9/8}}, \qquad t  \geq K.$$
Then $L := \int_K^\infty f_r(s) ds$ is a finite number and, as $t \to \infty$,
$$I[h A (H_c \cos)] - \int_0^K h A (H_c \cos) ds - \int_K^t f(s) ds - L = - \int_t^\infty f_r (s) ds = O(t^{-1/8}).$$
Thus, up to a constant, we can find the asymptotics of $I[h A (H_c \cos)]$ to $O(t^{-1/8})$ by simply integrating $f(t)$. This gives the result.
\proofend

The asymptotics of the first eight terms in (\ref{listoforder18}) to $O(t^{-1/16})$ now follow from lemma \ref{abelasymptoticslemma} in view of the asymptotics of lemmas \ref{HcHslemma} and \ref{lasteightlemma}. Indeed, all the functions acted on by $A$ in these terms are at most of $O(t^2)$ as $t \to \infty$ and lemma \ref{abelasymptoticslemma} applies.

It only remains to treat the last term in the expression (\ref{Ihhsexpression}) for $I_{hhs}$. This term is not in the form of an Abel integral. However, we can find its asymptotics using different arguments.

\begin{lemma}
The last term on the rhs of (\ref{Ihhsexpression}) satisfies
\begin{align}\label{lasttermI1234}
&2\int_0^t dt' h(t') \int_0^{t'}dt''  \frac{\cos{t''}H_c(t'') + \sin{t''}H_s(t'')}{\sqrt{ t- t' + t''}}
	\\ \nonumber
&\hspace{2cm} =  \sqrt{2\pi}\bigl[I_1(t)\cos{t} + I_2(t)\sin{t} + I_3(t)\cos{t} + I_4(t)\sin{t}\bigr] +O(t^{-1/2}),
\end{align}
where the functions $I_1, I_2, I_3, I_4$ are defined by
\begin{subequations}\label{I1234def}
\begin{align}
& I_1 = \int_0^tdt'' \frac{\cos{t''}}{\sqrt{t''}}\int_0^{t''} (H_c \cos^2 + (H_s + H_c)\cos \sin + H_s \sin^2) du,
	\\
& I_2 = \int_0^tdt'' \frac{\cos{t''}}{\sqrt{t''}}\int_0^{t''} (H_c \cos^2 + (H_s - H_c)\cos \sin - H_s \sin^2) du,
	\\
&I_3 = \int_0^tdt'' \frac{\sin{t''}}{\sqrt{t''}}\int_0^{t''} (-H_c \cos^2 + (-H_s + H_c)\cos \sin + H_s \sin^2) du,
	\\
&I_4 = \int_0^tdt'' \frac{\sin{t''}}{\sqrt{t''}}\int_0^{t''} (H_c \cos^2 + (H_s + H_c)\cos \sin + H_s \sin^2) du.
\end{align}
\end{subequations}
\end{lemma}
\proofbegin
For the purposes of finding the asymptotics to $O(t^{-1/2})$, the function $h(t')$ in the lhs of (\ref{lasttermI1234}) may be replaced with its leading asymptotics.
Indeed, lemma \ref{HcHslemma} implies that 
$$\cos{t''}H_c(t'') + \sin{t''}H_s(t'') = O(t''), \qquad t'' \to \infty.$$ 
Since
$$2\int_0^1 dt' \int_0^{t'}dt''  \frac{1}{\sqrt{ t- t' + t''}} = O(t^{-1/2}),$$
and
$$2\int_1^t dt' \frac{1}{t'^{7/2}}\int_0^{t'}dt''  \frac{t''}{\sqrt{ t- t' + t''}}  = O(t^{-1/2}),$$
it follows that to $O(t^{-1/2})$ the asymptotics of the lhs of (\ref{lasttermI1234}) coincide with those of the following expression:
\begin{align*}
&2\int_0^t dt' \sqrt{\frac{\pi }{2}} (\sin (t')+\cos (t')) \int_0^{t'}dt''  \frac{\cos{t''}H_c(t'') + \sin{t''}H_s(t'')}{\sqrt{ t- t' + t''}}
	\\
&- 2\int_1^t dt' \frac{1}{2 t'^{3/2}} \int_0^{t'}dt''  \frac{\cos{t''}H_c(t'') + \sin{t''}H_s(t'')}{\sqrt{ t- t' + t''}}.
\end{align*}
Changing the order of integration and performing the integral with respect to $dt'$, we can write this expression as
\begin{align} \nonumber
 2\pi \int_0^t dt'' & \bigl[\cos{t''}H_c(t'') + \sin{t''}H_s(t'')\bigr]
	\\\nonumber
& \times \biggl\{(\sin(t + t'') + \cos(t + t'')) \biggl(C\biggl(\sqrt{\frac{2}{\pi}}\sqrt{t}\biggr) - C\biggl(\sqrt{\frac{2}{\pi}}\sqrt{t''}\biggr)\biggr)
	\\\nonumber
& \qquad+ (\sin(t + t'') - \cos(t + t'')) \biggl(S\biggl(\sqrt{\frac{2}{\pi}}\sqrt{t}\biggr) - S\biggl(\sqrt{\frac{2}{\pi}}\sqrt{t''}\biggr)\biggr)\biggr\}
	\\\nonumber
+ 2 \int_0^{1}dt'' &\bigl[\cos{t''}H_c(t'') + \sin{t''}H_s(t'')\bigr] \frac{\sqrt{\frac{t''}{t}}-\sqrt{t+t''-1}}{t+t''}
	\\\label{lasttermmidstep}
- 2 \int_1^{t}dt'' &\bigl[\cos{t''}H_c(t'') + \sin{t''}H_s(t'')\bigr] \frac{t - t''}{(t + t'')\sqrt{t t''}}.
\end{align}
The first integral in (\ref{lasttermmidstep}) becomes after integration by parts
\begin{align*}
& \sqrt{2\pi}\int_0^t dt'' \biggl(\int_0^{t''} du [\cos{u}H_c(u) + \sin{u}H_s(u)](\sin(t + u) + \cos(t + u))\biggr)\frac{\cos{t''}}{\sqrt{t''}}
	\\
&+ \sqrt{2\pi}\int_0^t dt'' \biggl(\int_0^{t''} du [\cos{u}H_c(u) + \sin{u}H_s(u)](\sin(t + u) - \cos(t + u))\biggr)\frac{\sin{t''}}{\sqrt{t''}}.
\end{align*}	
After expanding the above trigonometric functions it follows that this expression equals 
$$ \sqrt{2\pi}\bigl[I_1(t)\cos{t} + I_2(t)\sin{t} + I_3(t)\cos{t} + I_4(t)\sin{t}\bigr].$$ 
It therefore only remains to prove that the last two integrals in (\ref{lasttermmidstep}) are of $O(t^{-1/2})$. For the next to last integral, this is clear. To see that the final integral is of $O(t^{-1/2})$, we use repeated integration by parts to write it as follows:
\begin{align} \nonumber
&- 2 \int_1^{t''} [\cos{u}H_c(u) + \sin{u}H_s(u)]du\frac{t - t''}{(t + t'')\sqrt{t t''}}\biggl|_{t''=1}^t
	\\\nonumber
&+ 2 \int_1^{t''} \int_1^{\tau} [\cos{u}H_c(u) + \sin{u}H_s(u)]dud\tau \frac{d}{dt''}\biggl(\frac{t - t''}{(t + t'')\sqrt{t t''}}\biggr)\biggl|_{t''=1}^t
	\\ \label{int1tppthreelines}
&- 2 \int_1^{t} dt'' \int_1^{t''} \int_1^{\tau} [\cos{u}H_c(u) + \sin{u}H_s(u)]dud\tau \frac{d^2}{dt''^2}\biggl(\frac{t - t''}{(t + t'')\sqrt{t t''}}\biggr).
\end{align}
The following facts imply that the last line of (\ref{int1tppthreelines}) is of $O(t^{-1/2})$:
\begin{itemize}
\item $\int_1^{t''} \int_1^{\tau} [\cos{u}H_c(u) + \sin{u}H_s(u)]dud\tau$ is of $O(t'')$ as $t'' \to \infty$,
\item $\frac{d^2}{dt''^2}\bigl(\frac{t - t''}{(t + t'')\sqrt{t t''}}\bigr) > 0$ for all $0 < t'' < t$,
\item $ \int_1^{t} dt'' t'' \frac{d^2}{dt''^2}\biggl(\frac{t - t''}{(t + t'')\sqrt{t t''}}\biggr) = O(t^{-1/2})$.
\end{itemize}
Consequently, the expression in (\ref{int1tppthreelines}) equals
\begin{align} \nonumber
- \frac{1}{t^2} \int_1^{t} \int_1^{\tau} [\cos{u}H_c(u) + \sin{u}H_s(u)]dud\tau + O(t^{-1/2}) = O(t^{-1/2}).
\end{align}
\proofend

We can now complete the proof of the theorem. Using the expressions for $H_c$ and $H_s$ of lemma \ref{HcHslemma}, the asymptotics of the integrals $\{I_j\}_1^4$ defined in (\ref{I1234def}) can be computed up to a constant. This yields the following formulas:
\begin{align}\nonumber
 & I_1 = \sqrt{\frac{\pi}{32}}\biggl\{t^{3/2} \sin {t}-\frac{\sqrt{t} \sin {t}}{2}-\frac{\sqrt{t} \sin {3t}}{6} +\frac{3\sqrt{t} \cos {t}}{2}  \biggr\} + k_9 + O(t^{-1/2}),
  	\\ \nonumber
 & I_2 = \sqrt{\frac{\pi}{32}}\biggl\{-\sqrt{t} \sin {t}-\frac{\sqrt{t} \cos {t}}{2} -\frac{\sqrt{t}\cos {3t}}{6} \biggr\}  + k_{10} + O(t^{-1/2}),
	\\ \label{I1234asym}
&  I_3 = \sqrt{\frac{\pi}{32}}\biggl\{ -\frac{\sqrt{t} \sin {t}}{2} +\frac{\sqrt{t} \sin {3t}}{6}-\sqrt{t}\cos {t}\biggr\}  + k_{11} + O(t^{-1/2}),
	\\ \nonumber
&  I_4 = \sqrt{\frac{\pi}{32}}\biggl\{-t^{3/2}\cos {t} +\frac{3\sqrt{t} \sin {t}}{2} -\frac{ \sqrt{t} \cos {t}}{2}+\frac{\sqrt{t} \cos {3t}}{6} \biggr\}  + k_{12} + O(t^{-1/2}),			
\end{align}
where $\{k_j\}_9^{12}$ are real constants.

Combining all of the above, we arrive at the following asymptotics of $T_j$, $j = 1, \dots, 7$, as $t \to \infty$:
\begin{align*}
T_1 & = \frac{1+i}{24}\left\{9(\cos{t} - \sin{t}) -\sqrt{3}(\cos{3t} - \sin{3 t})\right\}
+ O(t^{-1/2}),
	\\
T_2 &= \frac{1+i}{48}\left\{6 \sin {t}+ (\sqrt{3}-3) \sin {3t}-12 \cos{t}- (\sqrt{3}-3) \cos {3t} \right\} + O(t^{-1/2}),
	\\
T_3 & = -\frac{ \sin {t}}{4} (2 t - \sin {2t} - \cos {2t}) + O(t^{-3/2}),
	\\
T_4 &=\frac{1- i}{64\pi}\Bigl\{4 i \pi  t (\sin {t}-\cos {t}) 
	\\
& \hspace{1.7cm}+ 2 i \pi  \left(3 \sin {t}+\cos {t}+(\sqrt{3}-1) (\cos (3
   t)-\sin {3t})\right)  \Bigr\}  + O(t^{-5/2}),
   	\\
T_5 & = 	\frac{1- i}{64\pi}\Bigl\{-4 \pi  t (\cos {t}-3 \sin {t}) + 32 \sqrt{\tfrac{2}{\pi }} ((k_1+k_3) \cos{t}+(k_2-k_4) \sin {t})
	\\
& \hspace{1.7cm} +\sin {t} (\pi -32(k_5+k_6+k_7-k_8))+\cos {t} (\pi -32(k_5-k_6+k_7+k_8))
	\\
&\hspace{1.7cm} +2 \left(\sqrt{3}-1\right) \pi  (\cos {3t}-\sin {3t})\Bigr\}  + O(t^{-1/16}),
  	\\
T_6 &=	\frac{1- i}{64\pi}\Bigl\{ 4 i \pi  t (3 \sin {t}+\cos {t}) + 4 i \sin {t} (8 k_5+8 k_6-8 k_7+8 k_8-2 \pi )
	\\
&\hspace{.3cm} +4 i \cos {t} (8 k_5-8 (k_6+k_7+k_8))+2 i
   \left(\sqrt{3}-1\right) \pi  (\cos {3t}-\sin {3t})\Bigr\}  + O(t^{-1/16}),
   	\\
T_7 &=\frac{1- i}{64\pi}\Bigl\{  4 \pi  t (\sin {t}+\cos {t}) + \sin {t} (32 k_{10}+32 k_{12}-64 k_5-64 k_7-\pi)
	\\
&\hspace{1.7cm}+\cos {t} (32 k_{11}-64 k_5-64 k_7+32 k_9-\pi)
	\\
&\hspace{1.7cm} +2 \left(\sqrt{3}-1\right) \pi  (\cos {3t}-\sin {3t})\Bigr\} + O(t^{-1/16}).
\end{align*}

Substituting these expressions into (\ref{g13formula}), we find that $g_{13}$ satisfies (\ref{g11g12g13c}) with $c_1$ and $c_2$ given by
\begin{subequations}\label{c1c2def}
\begin{align}\nonumber
c_1 = \frac{1-i}{2\pi}\biggl\{&\sqrt{\frac{2}{\pi }} (k_1+k_3)-(3-i)k_5+(1-i) k_6-(3+i) k_7
   	\\ 
&  -(1+i) k_8+k_9 +k_{11}+\left(\frac{1}{8}+\frac{7 i}{16}\right) \pi\biggr\},
   	\\ \nonumber
c_2 = \frac{1-i}{2\pi}\biggl\{& \sqrt{\frac{2}{\pi }} (k_2-k_4)-(3-i) k_5-(1-i) k_6-(3+i) k_7
   	\\
&   +(1+i) k_8 + k_{10}+k_{12} -\left(\frac{1}{8}+\frac{11 i}{16}\right) \pi\biggr\},
\end{align}
\end{subequations}
where $\{k_j\}_1^{12}$ are defined in lemma \ref{lasteightlemma} and equation (\ref{I1234asym}).
This completes the proof of theorem \ref{g13th}. \proofend

\begin{remark}\upshape
1. Neither of the terms $T_3$-$T_7$ is asymptotically periodic by itself, however their sum is.

2. The leading asymptotic behavior of $O(t)$ of $T_5$, $T_6$, and $T_7$ can be obtained by simply replacing $h(t)$ by its leading asymptotics $\sqrt\frac{\pi}{2}(\cos{t} + \sin{t})$ in (\ref{Ixxxdef}) and computing all the integrals in terms of $C(z)$ and $S(z)$.
\end{remark}

\appendix
\section{An alternative system}\nequation
\renewcommand{\theequation}{A.\arabic{equation}}
Eliminating the unknown functions $g_1$ and $\bar{g}_1$ from the system (\ref{Goursat}), we find after simplification the following equations:
\begin{align*}
&  \beta M_2^2 + \biggl(L_2 + \frac{i\lambda}{2}\bar{g}_0 M_1\biggr)[M_{1t} - M_{1s} - 2g_0L_2] + \frac{i}{2}g_0M_2[M_{2t} + M_{2s} - 2\lambda \bar{g}_0 L_1]
\\
&  = M_2[L_{1t} - L_{1s}]
\end{align*}
and
\begin{align*}
 & \lambda \bar{\beta} M_1^2 + \biggl(L_1 - \frac{ig_0}{2}M_2\biggr)[M_{2t} + M_{2s} - 2\lambda \bar{g}_0 L_1]
  - \frac{i\lambda}{2}\bar{g}_0 M_1 [M_{1t} - M_{1s} - 2g_0 L_2] 
   	\\
&  = M_1[L_{2t} + L_{2s}].
\end{align*}
Replacing in the above $L_1$Ê and $L_2$ using (\ref{L1L2eqsa}) and (\ref{L1L2eqsb}), we find
\begin{align*}
 M_2(\partial_t - \partial_s) &\partial_s \int_{-t}^s \frac{M_1(t,\tau) d\tau}{\sqrt{s - \tau}}
- (M_{1t} - M_{1s})  \partial_s \int_{-t}^s \frac{M_2(t,\tau) d\tau}{\sqrt{s - \tau}}
 + \sqrt{\frac{\pi}{2}}e^{i\pi/4}g_0 (M_{2}^2)_s
	\\
&  + i\lambda|g_0|^2 \left(M_2 \partial_s \int_{-t}^s \frac{M_1(t,\tau) d\tau}{\sqrt{s - \tau}} - M_1 \partial_s \int_{-t}^s \frac{M_2(t,\tau) d\tau}{\sqrt{s - \tau}} \right)
	\\
& \qquad \qquad \qquad - \sqrt{\frac{2}{\pi}}e^{i\pi/4}g_0\left( \partial_s \int_{-t}^s \frac{M_2(t,\tau) d\tau}{\sqrt{s - \tau}}\right)^2 = 0
\end{align*}
and
\begin{align*}
 M_1(\partial_t + \partial_s)& \partial_s \int_{-t}^s \frac{M_2(t,s)d\tau}{\sqrt{s - \tau}}
 - (M_{2t} + M_{2s}) \partial_s \int_{-t}^s \frac{M_1(t,s)d\tau}{\sqrt{s - \tau}}
 + \lambda\sqrt{\frac{\pi}{2}}e^{i\pi/4}\bar{g}_0 (M_{1}^2)_s
	\\
&
+ i\lambda|g_0|^2 \left(M_2 \partial_s \int_{-t}^s \frac{M_1(t,s)d\tau}{\sqrt{s - \tau}} - M_1 \partial_s \int_{-t}^s \frac{M_2(t,s)d\tau}{\sqrt{s - \tau}}\right)
 	\\
& \qquad \qquad \qquad - \lambda\sqrt{\frac{2}{\pi}}e^{i\pi/4}\bar{g}_0 \left(\partial_s \int_{-t}^s \frac{M_1(t,s)d\tau}{\sqrt{s - \tau}}\right)^2
 = 0,
\end{align*}
This provides an alternative nonlinear system for $M_1$ and $M_2$. However, it is not clear if this system is perturbatively solvable. 


\section{Derivation of (\ref{Ffromcharacteristics1})}\label{characteristicapp}\nequation
\renewcommand{\theequation}{B.\arabic{equation}}
In order to solve equation (\ref{Fequation1}) we introduce the characteristic coordinates $\xi$ and $\eta$:
$$\begin{cases}
\xi = \frac{t+s}{2}, \\
\eta = \frac{t-s}{2} \end{cases} \qquad \Leftrightarrow \qquad \begin{cases} t = \xi + \eta, \\ s = \xi - \eta. \end{cases}$$
Thus,
$$\partial_t - \partial_s = \partial_\eta, \qquad \partial_t + \partial_s =\partial_\xi.$$
The restriction $-t < s < t$ implies that
$$0 < t-s < 2t, \qquad 0 < t+s < 2t,$$
thus
$$0 < \eta < T, \qquad 0 < \xi < T.$$
Equation (\ref{Fequation1}) yields
$$\frac{\partial F}{\partial \eta} =  f(\xi + \eta, \xi - \eta).$$
Hence,
$$F(t,s) = F(\xi + \eta_p, \xi - \eta_p) + \int_{\eta_p}^\eta f(\xi + \eta', \xi - \eta')d\eta', \qquad
0 \leq \eta_p \leq T.$$
Thus,
\begin{align}\label{Ftsetap}
F(t,s) = F\left(\frac{t+s}{2} + \eta_p , \frac{t+s}{2} - \eta_p\right) 
+ \int_{\eta_p}^{\frac{t-s}{2}} f\left(\frac{t+s}{2} + \eta', \frac{t+s}{2} - \eta'\right) d\eta'.
\end{align}
Letting $\eta_p = 0$ and employing the change of variables
$$\frac{t+s}{2} + \eta' = \tau,$$
equation (\ref{Ftsetap}) becomes equation (\ref{Ffromcharacteristics1}).

Alternatively, letting $\eta_p = 0$ and employing the change of variables
$$\frac{t+s}{2} - \eta' = \tau,$$
equation (\ref{Ftsetap}) yields
$$F(t,s) = F\left(\frac{t+s}{2}, \frac{t+s}{2}\right) - \int_{\frac{t+s}{2}}^s f(t+s - \tau, \tau)d\tau, \qquad
-t < s < t, \; 0 < t < T.$$

Similar considerations are valid for equation (\ref{Ffromcharacteristics2}).

\bigskip
\noindent
{\bf Note added:} After this paper was written we were pleased to find out about \cite{BSZ2003}, where the equation obtained from the KdV by adding the term $\alpha q(x,t)$, $\alpha > 0$, representing damping is considered. For this equation, it is shown via techniques of functional analysis, that if the small amplitude boundary forcing is periodic of period $T$, then the solution $q(x,t)$ eventually also becomes periodic of period $T$ at each spatial point. 

\bigskip
\noindent
{\bf Acknowledgement} {\it The authors acknowledge support from the EPSRC, UK. ASF acknowledges support from the Guggenheim foundation, USA.}

\bibliographystyle{plain}

\end{document}